\documentclass[11pt]{article}
\usepackage{amsmath, amsthm, amssymb}
\usepackage[dvips]{graphics,color}
\usepackage{epsfig}
\usepackage{latexsym}
\usepackage{graphicx}

\begin{document}

\title{A Little Statistical Mechanics for the Graph Theorist}

\author{\textit{Laura Beaudin}  \\ Saint Michael's College \\  Colchester, VT 05439 \\
\texttt{lbeaudin@smcvt.edu}
\and \textit{Joanna Ellis-Monaghan} \\ Saint Michael's College \\
Department of Mathematics \\ Colchester, VT 05439 \\ \texttt{jellis-monaghan@smcvt.edu}
\and \textit{Greta Pangborn} \\ Saint Michael's College \\ Department of Mathematics
\\ Colchester, VT 05439 \\ \texttt{gpangborn@smcvt.edu} \and
\textit{Robert Shrock} \\ Stony Brook University/SUNY \\ C. N. Yang Institute for Theoretical Physics
\\ Stony Brook, NY 11794 \\ \texttt{robert.shrock@sunysb.edu }}

\maketitle

\begin{abstract}

In this survey, we give a friendly introduction from a graph theory perspective to the $q$-state Potts model, an important statistical mechanics tool for analyzing complex systems in which nearest neighbor interactions determine the aggregate behavior of the system. We present the surprising equivalence of the Potts model partition function and one of the most renowned graph invariants, the Tutte polynomial, a relationship that has resulted in a remarkable synergy between the two fields of study.  We highlight some of these interconnections, such as computational complexity results that have alternated between the two fields.  The Potts model captures the effect of temperature on the system and plays an important role in the study of thermodynamic phase transitions.  We discuss the equivalence of the chromatic polynomial and the zero-temperature antiferromagnetic partition function, and how this has led to the study of the complex zeros of these functions.  We also briefly describe Monte Carlo simulations commonly used for Potts model analysis of complex systems. The Potts model has applications as widely varied as magnetism, tumor migration, foam behaviors, and social demographics, and we provide a sampling of these that also demonstrates some variations of the Potts model.  We conclude with some current areas of investigation that emphasize graph theoretic approaches.


\end{abstract}

\section{Introduction}

The Potts model of statistical mechanics models how micro-scale nearest neighbor energy interactions in a complex system determine the macro-scale behavior of the system.  This model plays an important role in the theory of phase transitions and critical phenomena in physics, and has applications as widely varied as adsorption of gases on substrates, tumor migration, foam behaviors, and social demographics. If we generalize the regular lattice (on which physicists normally consider the Potts model) to an abstract graph, the $q$-state Potts model partition function is an evaluation of one of the most renowned graph invariants, the Tutte polynomial [Tut47, 53, 67, 84].  The Potts model is in fact equivalent to the Tutte polynomial if both $q$ and the temperature are viewed as indeterminate variables.

Here, we give a friendly introduction to the interconnections between the $q$-state Potts model partition function of statistical mechanics and the Tutte and chromatic polynomials of graph theory.  In some respects, this paper complements the excellent survey of Welsh and Merino [WM00]. Where [WM00] is directed toward the physicist familiar with the Potts model who desires an introduction to the Tutte polynomial and its properties, here we hope to engage the graph theorist with an accessible introduction to the Potts model.  Ideally, this paper will generate further interest, particularly from the graph theoretical perspective that has proven so productive in this rapidly developing and important area.

The Potts model, from the mid-1900s, builds on the seminal work of Ernst Ising [Isi25].  The Ising model of magnetic behavior features nearest neighbor interactions between spins at each point on a lattice, where the spins can assume either of two values corresponding to magnetic polarization. From these local interactions, the aggregate global properties of the system can be studied. Central among these are phase transitions, that is, critical temperatures around which a small change in temperature results in an abrupt change in the magnetism of the system.  For the Ising model, there is no phase transition in the one dimensional case, and the transition has been determined exactly for the two dimensional square lattice.

Intrigued by a related model due to Ashkin and Teller [AT43], Cyril Domb suggested the study of what is now called the $q$-state Potts model to his Ph.D. student, Renfrey B. Potts, who developed the beginnings of the theory in his 1952 doctoral thesis [Pot52]. The Potts model generalizes the Ising model by allowing $q$ different spin values. Important thermodynamic functions such as internal energy, specific heat, entropy, and free energy, may be derived from the Potts model partition function. From a physics standpoint, one of the main reasons for the strong interest in the Potts model is that, for $q = 3$ and $q = 4$, it exhibits a continuous phase transition between high- and low-temperature phases with critical singularities in thermodynamic functions different from those of the Ising model. However, since its inception, myriad applications of the Potts model have emerged, and its usage now spans all areas of the sciences.

We present the essential concepts in this area, addressing the natural questions of: What is the Potts model and its partition function?  How are the Potts model partition function and the Tutte polynomial related? How does the Potts model capture phase transitions in thermodynamic functions?  What is the relationship between the Potts model and the chromatic polynomial? What is the computational complexity of the partition function?  How are Monte Carlo simulations used for the Potts model? Why is this model generating so much current interest? And finally, what are some current research directions in this field that emphasize graph theoretical approaches?

\section{The $q$-state Potts model}

Let $G$ be a graph and $S$ be a set of $q$ elements, called spins.  In the abstract, the spins may be numbers or colors, but typically they are values relevant to some specific application. For example, in studying uniaxial magnetic materials, $q = 2$, and the possible spins are $ + 1$ and $ - 1$.  In a foam model, there may be thousands of spins, one for each bubble in the foam.  In many applications, the graph $G$ is taken to be a regular lattice, but this assumption is not necessary.

A \emph{state} of a graph $G$ is an assignment of a single spin to each vertex of the graph. The \emph{Hamiltonian} is a measure of the energy of a state.  We begin with the simplest formulations, where the interaction energy (which may be thought of as a weight on each edge of the graph) is a constant $J$, and the Hamiltonian depends only on the nearest neighbor interactions (without any external field or other modifying forces).  The model is called \emph{ferromagnetic} if $J$ is positive and \emph{antiferromagnetic} if $J$ is negative. If $J$ is positive (respectively negative), the spin-spin interaction favors equal (respectively unequal) values of the spins on adjacent vertices.

We will see shortly that the two Hamiltonians in the following definition generate essentially equivalent formulations of the Potts model partition function.

\textbf{Definition 2.1.}Two common formulations of the Hamiltonian:
$$
h_1(\omega ) =  - J\sum\limits_{ij  \in E(G)} {\delta(\sigma _i ,\sigma _j)}
 \text{ and }  h_2(\omega) = J\sum\limits_{{ij} \in E(G)} {\left( {1 - \delta \left( {\sigma _i ,\sigma _j } \right)} \right)},
$$

where $\omega$ is a state of a graph $G$, $\sigma_i $ is the spin at vertex $i$, and $\delta $ is the Kronecker delta function.

\begin{figure}[htb]
  \centerline{\epsfig{file=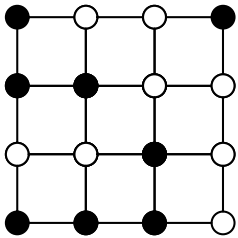, width=.25\textwidth}}
    \caption{A State of the $4 \times 4$ Square Lattice}
    \label{fig1}
 \end{figure}

For example, Figure~\ref{fig1} gives a graph state $\omega $ of the $4 \times 4$ square lattice with two choices of spin (black or white) at each vertex, with $h_1 (\omega ) =  - 11J$ and $h_2 (\omega ) = 13J$. Note that, up to the minus sign, $h_1 $ counts the edges with the same spins on their endpoints, and $h_2 $ counts the edges with different spins on their endpoints, so, for any state $\omega $, $h_2 \left( \omega  \right) = J\left| {E\left( G \right)} \right| + h_1 \left( \omega  \right)$.

\textbf{Definition 2.2.}  The Potts model partition function:  Given a set of $q$ spins and a Hamiltonian $h_i $ for $i = 1$  or $2$, the $q$-state Potts model partition function is $
P_i \left( G \right) = \sum {\exp ( - \beta (h_i (\omega )))} $.  Here the sum is over all possible states $\omega $ of $G$, and $\beta  = \kappa /T$, where $T$ represents the temperature of the system, and $\kappa  = 1.38 \times 10^{ - 23} $ joules/Kelvin is the Boltzmann constant.

The temperature $T = 1/(\kappa \beta )$ is an important variable in the model, although it need not represent physical temperature, but may be some other measure of agitation or volatility relevant to the particular application (economic factors in a sociological model for example).  Also, the product $J\beta $ occurs often and customarily appears in the physics literature as $J\beta  = K$; we will adopt this convention when convenient.

The Potts model partition function is the normalization factor for the Boltzmann probability distribution.  For systems such as the Potts model that follow Bolzmann distribution laws, the number of states with a given energy (Hamiltonian value) are exponentially distributed.  Thus, the probability of the system being in a particular state $\varpi $
 at temperature $T$ is:
 $Pr\left( {\varpi ,\beta } \right) = \exp ( - \beta h_i (\varpi ))/\sum {\exp ( - \beta h_i (\omega ))} $.

Since the two different formulations of the Potts model partition function are each natural to use in different contexts, the following observation that one is simply a scalar multiple of the other facilitates translating theoretical results from one context to the other.

\textbf{Observation 2.3.}

$\begin{array}{l}
 P_2 \left( {G;q,\beta } \right) = \sum {\exp \left( { - \beta h_2 \left( \omega  \right)} \right)}  =  \\
 \sum {\exp \left( { - \beta \left( {J\left| {E\left( G \right)} \right| + h_1 \left( \omega  \right)} \right)} \right)}  = \exp \left( { - K\left| {E\left( G \right)} \right|} \right)P_1 \left( {G;q,\beta } \right). \\
 \end{array}$

\begin{figure}[htb]
  \centerline{\epsfig{file=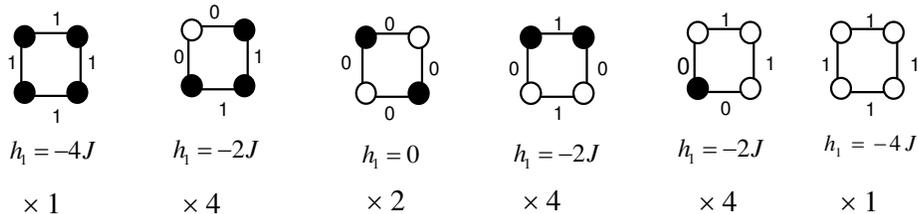, width=.8\textwidth}}
    \caption{Possible States for $C_4$}
    \label{fig2}
 \end{figure}

As a quick example, consider a single square, with two possible spins (white and black) at each vertex.  The possible states (up to rotation) and their Hamiltonians using $h_1$ are shown in Figure~\ref{fig2}.  Thus, the partition function for this graph is $P_1 \left( G \right) = 12\exp \left( {2K} \right) + 2\exp \left( {4K} \right) + 2$.

As noted following Definition 2.2, the probability of a particular state $\varpi $ occurring is $\exp ( - \beta h_i (\varpi ))/P_i (G)$.  Most importantly, since $\beta  = \kappa/ T$, this probability is temperature dependent.  Continuing the example in Figure~\ref{fig2}, we take $K = 1/T$, so $J = \kappa$  is positive, and consider the probability of the all-black state occurring at different temperatures.  Since $J$ is positive, and we are using $h_1 $, the all-black state is one of the two lowest energy states.  Thus, we would expect the system to favor the two low energy states equally at low temperatures, and be equally likely to be in any of the sixteen possible states at high temperatures.  If we let $\omega _b $ be the all-black state, then the probability of the system being in the all-black state, as a function of temperature, is $\Pr \left( {\omega _b ,T} \right) = \exp (4K)/(12\exp (2K) + 2\exp (4K)+ 2)$.  Evaluating this at various temperatures illustrates how the model captures the expected temperature dependent behavior: $\Pr (\omega _b ,10^{ - 2} ) = 1/2 = 0.50$, $\Pr (\omega _b ,2.29) = 0.19$; $\Pr (\omega _b ,10^5 ) = 1/16 = 0.0625$.

\section{Relating the Potts model and the Tutte polynomial}

From a graph theory perspective, one of the most remarkable aspects of the Potts model is its connection with the Tutte polynomial, one of the best known graph invariants.  See Fortuin and Kasteleyn [FK72] for the nascent stages of this discovery, and Wu [Wu82, Wu88] for further discussion from a physics viewpoint. More recent mathematical physics studies, such as Shrock [Shr00], Sokal [Sok00], and Welsh and Merino [WM00], give further exposition, while some relevant mathematical reviews include Tutte [Tut84], Biggs [Big93], Bollob\'{a}s [Bol98], and Welsh [Wel93].  The Tutte polynomial [Tut47, 53, 67] has a rich history and a wide range of applications.  We mention here only a very few properties necessary to establish the relation between it and the Potts model, and refer the reader to the relevant chapters of Welsh [Wel93] and Bollob\'{a}s [Bol98] and also to Brylawki [Bry80] and Brylawski and Oxley [BO92]  for an in-depth treatment of the Tutte polynomial, including generalizations to matroids.

The Tutte polynomial, $t(G;x,y)$, is a two-variable graph invariant that may be defined by a linear recursion relation in terms of deleting and contracting edges.  Recall that an edge is deleted from a graph $G$ by removing the edge, but not its incident vertices.  This is denoted by $G-e$.  An edge is contracted in $G$ by removing the edge and coalescing its incident vertices.  This is denoted by $G/e$.  An isthmus (a.k.a. a bridge, cut-edge, or co-loop) is an edge whose deletion increases the number of components of $G$.  A loop is an edge where both ends of the edge are incident with the same vertex.

\textbf{Definition 3.1.}  Deletion-contraction definition of the Tutte polynomial: $t(G;x,y) = t(G - e;x,y) + t( G/e;x,y)$ if $e$ is neither an isthmus nor a loop, and $
t\left( {G;x,y} \right) = x^i y^j $ if $G$ consists of $i$ isthmuses and $j$ loops.

\begin{figure}[htb]
  \centerline{\epsfig{file=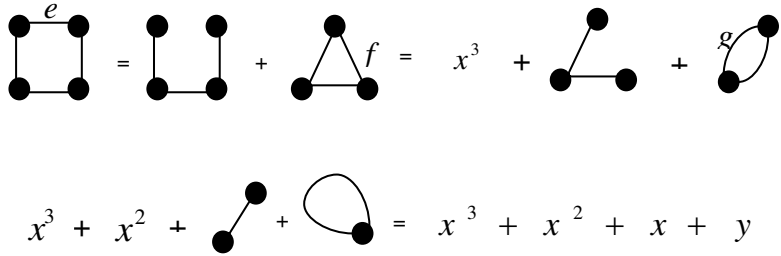, width=.7\textwidth}}
    \label{fig3}
    \caption{Calculation for $C_4$}
 \end{figure}

As a quick example, in Figure~\ref{fig3} we calculate that $t\left( {C_4 ;x,y} \right) = x^3  + x^2  + x + y$, where $C_4 $ is the cycle on four vertices (a square).  In the diagram, the labeled edge is deleted and contracted in the next step, and a graph consisting of only isthmuses and loops is evaluated as a monomial in $x$ and $y$ in the following step.

A surprising (and essential) property of the Tutte polynomial is that the result of the recursion process is independent of the order in which the edges are contracted and deleted.  Recall that for a graph $G = (V,E)$, a spanning subgraph is $G' = (V,A)$ where $A$ is a subset of $E$.  Also, $k(G)$ is the number of connected components in $G$ and $c(G)$ is the number of independent cycles in $G$.  One way to prove that $t(G;x,y)$ is independent of the order of contraction and deletion (and hence well-defined as given in Definition 3.1) is to express it as a sum of terms, each of which is evaluated for a given spanning subgraph.  This is analogous to the cluster representation for the Potts model in Fortuin and Kasteleyn [FK72].  In particular, $t(G;x,y) = \sum\limits_{G'} {(x - 1)^{k(G') - k(G)} (y - 1)^{c(G')} } $, where the sum is over the spanning subgraphs $G'$ of $G$ (see Bollob\'{a}s [Bol98], and Biggs [Big93]). Since the set of spanning subgraphs is unique, this proves that $t(G;x,y)$ is well defined for a given graph $G$.  Alternatively, using induction on the number of edges, one can show that $t(G;x,y)$ is equivalent, up to a prefactor, to the dichromatic polynomial, defined by Tutte [Tut67,47,53] as $Z\left( {G;u,v} \right) = \sum\limits_{A \subseteq E\left( G \right)}^{} {u^{k(A)} v^{\left| A \right|-\left| V(G) \right|+ k(A) } } $.  Specifically,
\begin{equation}
u^{k(G)} v^{\left| V \right| - k\left( G \right)} t\left( {G;\frac{{u + v}}{v},v + 1} \right) = Z\left( {G;u,v} \right).
\end{equation}

Theorem 3.1 gives another critical property of the Tutte polynomial.  The Tutte polynomial is universal in that essentially any other multiplicative graph invariant that has a deletion-contraction reduction must be an evaluation of it.

\textit{\textbf{Theorem 3.1. }  If $f(G)$ is a function on graphs with $a, b$ having $ab \neq 0$, such that
\begin{description}
\item[A.]   $f(G) = 1$ if $G$ consists of only one vertex and no edges,
\item[B.]   $f(G) = af(G - e) + bf(G/e)$ whenever e is neither a loop nor a bridge,
\item[C.]   $f(GH) = f(G)f(H)$ where either $GH $ is the disjoint union of $G$ and $H$ or $G$ and $H$ share at most one vertex,
\end{description}}

\textit{then $f$ is an evaluation of the Tutte polynomial of the form
$f(G) = a^n b^r t(G;x_0/b,y_0/a)$.
Here $n$ is the nullity of $G$, that is, $n = \left| {E(G)} \right| - \left| {V(G)} \right| + k(G)$, and $r$ is the rank of $G$, that is, $r = \left| {V(G)} \right| - k(G)$.  Also, $x_0  = f(K_2 )$ and $y_0  = f(L)$, where $K_2 $
 is the complete graph on two vertices, and $L$ is a single vertex with a single loop edge.}

This universality provides the basis for a proof of the connection between the Tutte polynomial and the Potts model partition function.  If we consider the $h_1 $ Hamiltonian in the Potts model intuitively, we note that, in a given state of the graph, if the endpoints of an edge have different spins, then the Kronecker delta value is zero, and the edge contributes nothing, so it might as well be deleted.  On the other hand, if the endpoints have the same spin, they interact with the neighboring points in exactly the same way, so they may be coalesced, with the edge contracted.  However, this edge does contribute to the Hamiltonian, so there is a weighting factor when the edge is contracted.

Although the Potts model partition function seems likely to satisfy the conditions of Theorem 3.1, and thus be an evaluation of the Tutte polynomial, it does not satisfy condition A: if there are $q$ spins, the Potts model partition function of a single vertex is $q$, not 1.  However, there is a very common device for applying Theorem 3.1 in such situations, namely introducing a factor of some term raised to the power of $k\left( G \right)$, a trick employed in the proof of the following theorem.  Although the connection between the Potts model partition function and the Tutte polynomial was first recognized by Fortuin and Kastelyn [FK72], the proof below is modeled on those in Welsh [Wel93] and Bollob\'{a}s [Bol98].


\textit{\textbf{Theorem 3.2}
Let $\tilde P\left( {G;q,\beta } \right) = q^{ - k\left( G \right)} P_1 \left( {G;q,\beta } \right)$.}

\textit{Then $\tilde P\left( {G;q,\beta } \right) = v^{\left| {V(G)} \right| - k(G)} t\left( {G;(q + v)/v,v + 1} \right)$,} 

\textit{and thus $P_1 \left( {G;q,\beta } \right) = q^{k(G)} v^{\left| {V(G)} \right| - k(G)} t\left( {G;(q + v)/v,v + 1} \right)$, where $v = \exp (J\beta ) - 1$. }

\textit{Proof:} For ease of reading, we will suppress the subscript 1 on $P$ and $h$ in the following proof.  The proof consists simply of verifying that $\tilde P$ satisfies the conditions of Theorem 3.1.  $\tilde P$ clearly satisfies condition A.

For condition B, let $e = \left\{ {c,d} \right\}$
 be an edge of $G$ which is neither a loop nor an isthmus, and write $s\left( c \right)$
 and $s\left( d \right)$
 for the spins at c and d respectively.  Then
$$
 \tilde P\left( {G;q,\beta } \right) = q^{ - k\left( G \right)} \sum\limits_{\omega  \in {\text{ states of }}G} {\exp \left( { - \beta h\left( \omega  \right)} \right)}
$$

\[
= q^{ - k\left( G \right)} \left( {\sum\limits_{\scriptstyle \omega  \in {\text{ states of }}G \hfill \atop
  \scriptstyle {\text{with }}s\left( c \right) \ne s\left( d \right) \hfill} {\exp \left( { - \beta h\left( \omega  \right)} \right)} } \right)
\tag{2}  + q^{ - k\left( G \right)} \left( {\sum\limits_{\scriptstyle \omega  \in {\text{ states of }}G \hfill \atop
  \scriptstyle {\text{with }}s\left( c \right) = s\left( d \right) \hfill} {\exp \left( { - \beta h\left( \omega  \right)} \right)} } \right).
\]


Note that if $\omega _G $ is any state of $G$, then there is a unique state $\omega _{G - e} $
 of $G - e$ where each vertex has the same spin as it has in $\omega _G $.  Also, if $\omega _G $
is any state of $G - e$ with $s\left( c \right) = s\left( d \right)$, then there is a unique state $\omega _{G/e} $ of $G/e$, where the vertex resulting from identifying $c$ and $d$ has the common value $s\left( c \right) = s\left( d \right)$, and each other vertex has the same spin as it does in $\omega _{G - e} $.  Furthermore, if $s\left( c \right) \ne s\left( d \right)$, then $h\left( {\omega _G } \right) = h\left( {\omega _{G - e} } \right)$; and if $s\left( c \right) = s\left( d \right)$, then $h\left( {\omega _G } \right) = h\left( {\omega _{G - e} } \right) + J$.  Thus, equation (2) becomes

$$
q^{ - k\left( G \right)} \left( {\sum\limits_{\scriptstyle \omega  \in {\text{ states of }}G - e \hfill \atop
  \scriptstyle {\text{with }}s\left( c \right) \ne s\left( d \right) \hfill} {\exp \left( { - \beta h\left( \omega  \right)} \right)} } \right){\text{     }}
\\+ q^{ - k\left( G \right)} \exp \left( {J\beta } \right)\left( {\sum\limits_{\scriptstyle \omega  \in {\text{ states of }}G - e \hfill \atop
  \scriptstyle {\text{with }}s\left( c \right) = s\left( d \right) \hfill} {\exp \left( { - \beta h\left( \omega  \right)} \right)} } \right).
$$

The right-most term is nearly $\tilde P\left( {G - e;q,\beta } \right)$, since $e$ being neither a bridge nor a loop means that $k\left( G \right) = k\left( {G - e} \right)$, but we are missing the states of $G - e$ where $s\left( c \right) = s\left( d \right)$.  So we simply add and subtract them, getting
$$
 q^{ - k\left( {G - e} \right)} \left( {\sum\limits_{\scriptstyle \omega  \in {\text{ states of }}G - e \hfill \atop
  \scriptstyle {\text{with }}s\left( c \right) \ne s\left( d \right) \hfill} {\exp \left( { - \beta h\left( \omega  \right)} \right)} } \right)
+ q^{ - k\left( {G - e} \right)} \exp \left( {J\beta } \right)\left( {\sum\limits_{\scriptstyle \omega  \in {\text{ states of }}G - e \hfill \atop
  \scriptstyle {\text{with }}s\left( c \right) = s\left( d \right) \hfill} {\exp \left( { - \beta h\left( \omega  \right)} \right)} } \right)
$$
$$
  - q^{ - k\left( {G - e} \right)} \left( {\sum\limits_{\scriptstyle \omega  \in {\text{ states of }}G - e \hfill \atop
  \scriptstyle {\text{with }}s\left( c \right) = s\left( d \right) \hfill} {\exp \left( { - \beta h\left( \omega  \right)} \right)} } \right)
+ q^{ - k\left( {G - e} \right)} \left( {\sum\limits_{\scriptstyle \omega  \in {\text{ states of }}G - e \hfill \atop
  \scriptstyle {\text{with }}s\left( c \right) = s\left( d \right) \hfill} {\exp \left( { - \beta h\left( \omega  \right)} \right)} } \right)
$$
$$
 = q^{ - k\left( {G - e} \right)} \left( {\exp \left( {J\beta } \right) - 1} \right)\left( {\sum\limits_{\scriptstyle \omega  \in {\text{ states of }}G - e \hfill \atop
  \scriptstyle {\text{with }}s\left( c \right) = \left( d \right) \hfill} {\exp \left( { - \beta h\left( \omega  \right)} \right)} } \right)
 + q^{ - k\left( {G - e} \right)} \left( {\sum\limits_{\omega  \in {\text{ states of }}G - e} {\exp \left( { - \beta h\left( \omega  \right)} \right)} } \right). \\
$$

The second term is now $\tilde P\left( {G - e;q,\beta } \right)$.  For the first term, note that since $e$ is neither a bridge nor a loop, $k\left( {G - e} \right) = k\left( {G/e} \right)$.  Also, the states of $G - e$ with $s\left( c \right) = s\left( d \right)$ correspond exactly to the states of $G/e$, and furthermore a state of $G - e$ with $s\left( c \right) = s\left( d \right)$ has the same Hamiltonian as the corresponding state of $G/e$.  Thus, the first term becomes
$$
q^{ - k\left( {G/e} \right)} \left( {\exp \left( K \right) - 1} \right)\sum\limits_{\omega  \in {\text{ states of }}G/e} {\exp \left( { - \beta h\left( \omega  \right)} \right)}  = \left( {\exp \left( K \right) - 1} \right)\tilde P\left( {G/e;q,\beta } \right).
$$

Thus, if $e$ is neither a bridge nor a loop, then $\tilde P\left( {G;q,\beta } \right) =$ $ \tilde P\left( {G - e;q,\beta } \right) +$ $ \left( {\exp (K) - 1} \right)$ $\tilde P\left( {G/e;q,\beta } \right)$, which satisfies condition B with $a = 1$ and $b = \exp \left( K \right) - 1$.

For condition C we write $G \cup H$ when $G$ and $H$ are disjoint, and $G * H$ when $G$ and $H$ share a single vertex.  Condition C is easily satisfied when $G$ and $H$ are disjoint, since in this case a state $\omega _{G \cup H} $ of $G \cup H$ is just two independent states $\omega _G $ and $\omega _H $ of $G$ and H respectively, and $h\left( {\omega _{G \cup H} } \right) = h\left( {\omega _G } \right) + h\left( {\omega _H } \right)$.  With this, and noting that $k\left( {G \cup H} \right) = k(G) + k(H)$, it follows that
$$
\tilde P\left( {G \cup H;q,\beta } \right) = q^{ - k\left( {G \cup H} \right)} \sum\limits_{\scriptstyle {\text{states }}\omega  \hfill \atop
  \scriptstyle {\text{of }}G \cup H \hfill} {\exp \left( { - \beta h_1 \left( \omega  \right)} \right)}  =
$$
$$
q^{ - k\left( G \right) - k\left( H \right)} \sum\limits_{\scriptstyle {\text{states }}\omega _G {\text{ of }}G{\text{ and}} \hfill \atop
  \scriptstyle {\text{states }}\omega _H {\text{ of }}H \hfill} {\exp \left( { - \beta \left( {h\left( {\omega _G } \right) + h\left( {\omega _H } \right)} \right)} \right)}  =
$$
$$
q^{ - k\left( G \right)} \sum\limits_{{\text{states }}\omega _G {\text{ of }}G{\text{ }}} {\exp \left( { - \beta h\left( {\omega _G } \right)} \right)} {\text{   }}q^{ - k\left( H \right)} \sum\limits_{{\text{states }}\omega _H {\text{ of }}H} {\exp \left( { - \beta h\left( {\omega _H } \right)} \right) = }
$$
$$
\tilde P\left( {G;q,\beta } \right)\tilde P\left( {H;q,\beta } \right).
$$

In the case where $G$ and $H$ share a single vertex $u$, a state $\omega _{G * H}$ of $G * H$ corresponds to two states $\omega _G $ and $\omega _H $ which both have the same spin at $u$.  Here, $k\left( {G * H} \right) = k(G) + k(H) - 1$.  Note that the number of states of $H $ in which $s(u)=a$  is equal to the number of states in which $s(u)=b$, for any other spin $b$, since we can simply exchange the roles of $a$ and $b$ in any state.  Thus,  

$$\sum\limits_{\omega _H  \in {\text{ states of }}H} {\exp \left( { - \beta h\left( {\omega _H } \right)} \right)}  = q\sum\limits_{\scriptstyle \omega _H  \in {\text{ states of }}H \hfill \atop
  \scriptstyle {\text{with }}s(u) = a \hfill} {\exp \left( { - \beta h\left( {\omega _H } \right)} \right)}.$$    

With this we have:

$$
\tilde P\left( {G * H;q,\beta } \right) = q^{ - k\left( {G * H} \right)} \sum\limits_{\scriptstyle {\text{states }}\omega  \hfill \atop
  \scriptstyle {\text{of }}G * H \hfill} {\exp \left( { - \beta h\left( \omega  \right)} \right)}  =
$$

$$
q^{ - k\left( G \right) - k\left( H \right) + 1} \sum\limits_{\scriptstyle \omega _H  \in {\text{ states of }}H \hfill \atop
  {\scriptstyle {\text{with }}s\left( u \right){\text{ in }}\omega _H  \hfill \atop
  \scriptstyle {\text{equal to }}s\left( u \right){\text{ in }}\omega _G  \hfill}} {\exp \left( { - \beta h\left( {\omega _G } \right)} \right)} \exp \left( { - \beta h\left( {\omega _H } \right)} \right) =
$$
$$
q^{ - k\left( G \right) - k\left( H \right) + 1} \sum\limits_{\omega _G  \in {\text{ states of }}G} {\left( {\sum\limits_{\scriptstyle \omega _H  \in {\text{ states of }}H \hfill \atop
  {\scriptstyle {\text{with }}s\left( u \right){\text{ in }}\omega _H  \hfill \atop
  \scriptstyle {\text{equal to }}s\left( u \right){\text{ in }}\omega _G  \hfill}} {\exp \left( { - \beta h\left( {\omega _G } \right)} \right)\exp \left( { - \beta h\left( {\omega _H } \right)} \right)} } \right)}  =
$$
$$
q^{ - k\left( G \right) - k\left( H \right) + 1} \sum\limits_{\omega _G  \in {\text{ states of }}G} {q^{ - 1} \sum\limits_{\omega _H  \in {\text{ states of }}H} {\exp \left( { - \beta h\left( {\omega _G } \right)} \right)\exp \left( { - \beta h\left( {\omega _H } \right)} \right)} }  =
$$
$$q^{ - k\left( G \right) - k\left( H \right)} \sum\limits_{\omega _G  \in {\text{ states of }}G} {\exp \left( { - \beta h\left( {\omega _G } \right)} \right)} \sum\limits_{\omega _H  \in {\text{ states of }}H} {\exp \left( { - \beta h\left( {\omega _H } \right)} \right)}
 =
$$
$$
 \tilde P\left( {G;q,\beta } \right)\tilde P\left( {H;q,\beta } \right).
 $$

Thus, condition C is satisfied, and Theorem 3.1 applies to $\tilde P$.  To apply the conclusion of Theorem 3.1, it remains to compute  $\tilde P$ at a single isthmus and a single loop, that is, find $\tilde P\left( {K_2 } \right)$ and $\tilde P\left( L \right)$.

For a loop, note that there are $q$ states, and since both endpoints of a loop necessarily have the same value, the Hamiltonian of each state is 1.  Thus,
$$
\tilde P\left( {L;q,\beta } \right) = q^{ - 1} \sum\limits_{q{\text{ states}}} {\exp \left( { - \beta \left( { - J \cdot 1} \right)} \right)}  = q^{ - 1} q\exp \left( K \right) = \exp \left( K \right).
$$

For $K_2 $, there are $q$ states where the spins on the endpoints are equal, giving a Hamiltonian of 1.  Then there are $q\left( {q - 1} \right)$ states where the spins on the endpoints are different, giving a Hamiltonian of 0.  Thus, $$
\tilde P\left( {B;q,\beta } \right) = q^{ - 1} \left( {q\left( {q - 1} \right)\exp \left( { - K \cdot 0} \right) + q\exp \left( {K \cdot 1} \right)} \right) = \left( {\exp \left( K \right) + q - 1} \right).
$$

We now apply Theorem 3.1 with $a = 1$, $b = \left( {\exp \left( K \right) - 1} \right)$, $x_0  = \left( {\exp \left( K \right) + q - 1} \right)$, and
$y_0  = \exp \left( K \right)$.  To simplify the expression, we first set $\nu  = \exp \left( K \right) - 1$, so $a = 1$, $b = v$, $y_0  = v + 1$, and $x_0  = q + v$.  Now Theorem 3.1 yields
$$
\tilde P\left( {G;q,\beta } \right) = v^{\left| {V(G)} \right| - k(G)} t\left( {G;\frac{{q + v}}{v},v + 1} \right).
$$

Thus, since $\tilde P\left( {G;q,\beta } \right) = q^{ - k\left( G \right)} P_1 \left( {G;q,\beta } \right)$, it follows that
$$
P_1 \left( {G;q,\beta } \right) = q^{k(G)} v^{\left| {V(G)} \right| - k(G)} t\left( {G;\frac{{q + v}}{v},v + 1} \right).   \text{       ///}
$$

As an example, we recall from Figure~\ref{fig2} that $P_1 \left( G \right) $ $= 12\exp (2K) + 2\exp (4K) + 2
$, so $P_1 \left( {K_2 ,2,\beta } \right) = 12\left( {v + 1} \right)^2  + 2\left( {v + 1} \right)^4  + 2$, since there were $q = 2$ spins in the example.  Also recall from Figure~\ref{fig3} that $t(C_4;x,y)=x^3+x^2+x+y$.  A quick calculation verifies that  $2^1 v^{4 - 1} t\left( {G;\frac{{2 + v}}{v},v + 1} \right) = 12\left( {v + 1} \right)^2  + 2\left( {v + 1} \right)^4  + 2$.

We then have the following two immediate corollaries, with Corollary 3.3 following from Observation 2.3, and Corollary 3.4 from equation (1).

\textbf{Corollary 3.3.}
$$
P_2 \left( {G;q,\beta } \right) = q^{k\left( q \right)} \left( {v + 1} \right)^{ - \left| {E\left( G \right)} \right|} v^{\left| {V\left( G \right)} \right| - k\left( G \right)} t\left( {G;\frac{{q + v}}{v},\exp \left( K \right)} \right).
$$

\textbf{Corollary 3.4.}
$$
P_1 \left( {G;q,\beta } \right) = Z\left( {G;q,v} \right) = \sum\limits_{A \subseteq E\left( G \right)}^{} {q^{k\left( A \right)} v^{\left| A \right|} }.
$$

Thus, if $q$ and $v$ are viewed as indeterminates, the Potts model partition function is exactly equal to the dichromatic polynomial $Z\left( {G;q,v} \right)$, and in fact the partition function is typically denoted by $Z$ in the physics literature.  Corollary 3.4  also leads to the property (not apparent from the original definition) that the Potts model partition function is a two-variable polynomial with maximal degree in $q$ equal to the number of vertices of $G$ and maximal degree in $v$ equal to the number of edges of $G$.  Theorem 3.2 and Corollary 3.4 also show how the Tutte polynomial may be thought of as an analytic continuation of the Potts model with its positive integer values for $q$.  The random cluster model of Fortuin and Kasteleyn [FK72] also extends the Potts model in this way to $\mathbb{R}^+ $, and is likewise an evaluation of the Tutte polynomial.

One common extension of the Potts model, dating to the earliest work in the quantum theory of magnetism in the 1920's and 1930's, involves allowing the interaction energy to depend on the specific edge, rather than be constant throughout.  With this, the Hamiltonian becomes $h(\omega ) = \sum\limits_{\{ ij\}  \in E(G)} {J_{ij} \delta (\sigma _i ,\sigma _j )} $, where $J_{ij} $ (or $J_e $)  is the interaction energy on the edge $e=\{i,j\}$.  The partition function then becomes $P(G) = Z(G;q,{\bf{v}}) = \sum\limits_{A \subseteq E(G)} {q^{k\left( A \right)} \prod\limits_{e \in A} {v_e } } $, where  $v_e=\textnormal{exp}(\beta J_e)-1$.  See, for example, Fortuin and Kasteleyn [FK72] and Baxter [Bax82] for an edge weight generalization of the Potts model, and also recent work by Sokal [Sok00, 01a].

In recent years, the Tutte polynomial has also been extended to incorporate edge weights.  Here, these weights may also depend on whether the edge is contracted, deleted, or evaluated as an isthmus or a loop as the polynomial is recursively computed.  In the most general case, however, care must be taken with a set of relations on small graphs (two edges in parallel, a cycle with three edges, and three edges in parallel) to assure that the resulting function is well defined. See Traldi [Tra89], Zaslavsky [Zas92], Bollob\'{a}s and Riordan [BR99], and Ellis-Monaghan and Traldi [E-MT06].  The generalized partition function given above satisfies these relations, however, and thus the connection between the Potts model partition function and the Tutte polynomial extends to systems with edge-dependent interaction energies.

\section{Thermodynamic functions and phase transitions}

An important goal of statistical mechanics is to determine phase transitions, that is, critical temperatures around which a small change in temperature results in an abrupt, nonanalytic change in various physical properties. Roughly speaking, this phase transition temperature separates the two phases of the system.  For temperatures above this critical temperature, the system (in the absence of an external biasing field) exhibits no long-range order, e.g., no spontaneous magnetization in the case of a magnetic system, while for temperatures below the critical temperature it does exhibit such order. (See Stanley [Sta71], and Plischke and Bergesen [PB06].)
In the following discussion, we elide technical caveats concerning such details as the choice and growth of the lattices, interchanging limit signs, boundary conditions, and convergence, etc., in order to provide a broad picture of the general principles.  However, treating these technicalities carefully can present significant challenges in determining phase transitions for various applications.  Bearing this in mind, when we speak of taking the thermodynamic (or infinite volume) limit below, we mean specifying an appropriate infinite family of graphs, such as square lattices,  and taking the limit of an expression as the size of the graphs goes to infinity.

In the Potts model, important thermodynamic functions such as internal energy, specific heat, entropy, and free energy (denoted $U$, $C$, $S$, and $F$, respectively) may all be derived from the partition function, $Z = Z\left( {G;q,v} \right)$. For example, the internal energy, which is the sum of the potential and kinetic energy and is defined as $U = \frac{1}{Z}\sum {h\left( \omega  \right)\exp \left( { - \beta h\left( \omega  \right)} \right)} $, may be expressed by $\frac{{ - \partial \ln \left( Z \right)}}{{\partial \beta }}$.  Specific heat, or the energy required to raise a unit amount of material one temperature increment, is $C = \frac{{\partial U}}{{\partial T}}$.  Entropy, a measure of the randomness and disorder in a system, is $S =  - \kappa \beta \frac{{\partial \ln \left( Z \right)}}{{\partial \beta }} + \kappa \ln \left( Z \right)$.  Finally, the total free energy is $F = U - TS =  - \kappa T\ln \left( Z \right)$.

It is convenient to work with the dimensionless, reduced free energy $f =  - \beta F$, so the reduced free energy per unit volume (or in this context, per vertex) is
$$
f\left( {G;q,\beta } \right) = \frac{1}{{\left| {V\left( G \right)} \right|}}\ln \left( {P_1 \left( {G;q,\beta } \right)} \right).
$$

For a fixed graph $G$, this is clearly an analytic function in both $q$ and $T$.  Any failures of analyticity can only occur in the infinite volume limit, that is
$$
f\left( {\Gamma ;q,\beta } \right) = \mathop {\lim }\limits_{n \to \infty } \frac{1}{{\left| {V\left( {G_n } \right)} \right|}}\ln \left( {P_1 \left( {G_n ;q,\beta } \right)} \right),
$$
where $\Gamma $ is an appropriate infinite family of graphs.

If the resulting limit of the reduced free energy is an analytic function of $T$ (generally $q$ is fixed), then the model has no phase transitions since $U$, $C$, $F$, and $S$ will be analytic as well.  This happens for example for finite temperature $\beta  \in [0,\infty )$ in the Ising model on the one-dimensional lattice.  Values of $T$ where analyticity fails are critical temperatures corresponding to phase transitions.  Since phase transitions are manifest as failures of analyticity in the thermodynamic limit of the reduced free energy, the goal is either to determine these points of nonanalyticity or to establish analyticity in some region, and here the behavior of the partition function is the key. Furthermore, if we want to know how the thermodynamic functions behave near a critical temperature, again understanding the partition function is essential.

If $T_c $ is a critical temperature, and we write $\tau  = (T - T_c )/T_c$, then the goal is to express any of these thermodynamic functions in the form $g\left( \tau  \right)$, where $g\left( \tau  \right)$ is roughly equal to $c\left| \tau  \right|^p $ near $\tau  = 0$ for some $p$ called the critical exponent.  The critical exponents fall into a set of discrete universality classes, where universality means that the values of the critical exponents are independent of parameters such as the interaction energy $J$ and the choice of lattice (although not its dimension).

Phase transitions are broadly classified as first-order (discontinuous) or second-order (continuous, but nonanalytic). Second-order phase transitions are further classified by the critical exponents, since if the phase transitions are continuous, the leading singular behavior as $\tau  \to 0$ of the thermodynamic quantities may normally be written in the form $c\left| \tau  \right|^p $, where $p$ is a positive or negative power (there are also cases where the singularity is non-algebraic).  See Stanley [Sta71], Fisher [Fis74], and Plischke and Bergesen [PB06].

One of the important features of the two-dimensional ferromagnetic $q$-state Potts model is that, for the thermodynamic limit of regular two-dimensional lattice graphs, it provides, within one model, a set of several different universality classes associated with second-order phase transitions depending on $q$, in particular for $q = 3$ and $q = 4$, which generalize the Ising $q = 2$ case.  For $q \ge 5$, the phase transition of this two-dimensional ferromagnetic Potts model is first-order. Onsager [Ons44] calculated an exact closed-form expression for the free energy of the Ising model on a square lattice (in the absence of an external magnetic field) in 1944  (reviewed by McCoy and Wu in [MW73]).  For values of $q \ge 3$
, the free energy of the $q$-state Potts model for arbitrary temperatures (either with or without an external magnetic field) has never been calculated exactly for (the infinite volume limit of)  a lattice of dimension two or more.

Just as studying the generalization to complex variables can lead to greater insight into properties of functions of real variables, the generalization of the variable $K = J/(\kappa T)$ from real values (positive for ferromagnetic, negative for antiferromagnetic) to complex values has proven quite informative.  Indeed, this is necessary in order to study zeros of the partition function. For fixed $q$, the accumulation set of the zeros of the partition function in the infinite volume limit form curves in the complex plane referred to as the phase diagram. For example, Fisher [Fis65] showed that for the Ising model on the square lattice, as the number of vertices goes to infinity in the infinite volume limit, the zeros of the partition function asymptotically merge to form two circles in the complex plane, $|v| = \sqrt {2}$ and $|v + 2| = \sqrt {2}$. (Earlier, Lee and Yang had studied the corresponding extension to complex numbers of the magnetic field [YL52, LY52].)  This leads to the study of regions in the plane of the complex-temperature variable $v$ which are analytic continuations of the physical high- and low-temperature phases of the ferromagnetic and antiferromagnetic models.  Determinations of the complex-temperature phase diagram for the $q$-state Potts model for arbitrary $q$ were given for finite-width, infinite-length strips in Shrock [Shr00], Chang and Shrock [CS00, 01b], Chang, Salas and Shrock [CSS02], Chang, Jacobsen, Salas and Shrock [CJSS04], Jacobsen, Richard and Salas [JRS06], and the shape of these phase diagrams in general is discussed by Biggs [Big02b].

The determination of these phase diagrams requires a particularly interesting combination of methods from mathematical physics, graph theory, complex analysis, and algebraic geometry (since the phase boundaries are algebraic curves).  Given the exact partition function as a function of both arbitrary $q$ and $v$, it is possible to determine regions of analyticity in the complex $q$ plane as a function of the temperature variable $v$ for both the ferromagnetic and antiferromagnetic cases.  In addition to determining the phase diagram in $v$ for fixed $q$ and in $q$ for fixed $v$, one can also determine it when $q$ and $v$ satisfy a given functional relation (see Chang and Shrock  [CS06]). Complex-temperature zeros of Potts models for fixed values of $q$ (beyond the known Ising case) have also been calculated on finite patches of two-dimensional lattices (e.g. Martin and Maillard [MM86], Martin [Mar91], Chen, Hu and Wu [CHW96], Matveev and Shrock [MS96], Kim and Creswick [KC01]).

Additionally, Salas and Sokal [SS97] have shown that the antiferromagnetic Potts model with $q$ spins has no finite temperature phase transition on lattices where each vertex has degree less than q/2.

Quite a lot is known in the ferromagnetic case as opposed to the antiferromagnetic case. For example, the value of the critical temperature for the ferromagnetic Potts model on the infinte volume limit of the square lattice has been determined to be $\kappa T_c  = J/\ln (1 + \sqrt q )$, and the critical behavior of the two-dimensional Potts ferromagnetic model is known (see Baxter [Bax82] and Wu [Wu82]). Simulations on the square lattice agree with this formula for the critical temperature, and there have been a number of studies of the Potts model on various 2- and 3-dimensional lattices, giving valuable approximations with sufficient accuracy for relevant applications. Wu [Wu82, 84] and Salas and Sokal [SS97], for example, provide a survey of results and approximations.   Further insight into the critical exponents has been gained from the use of conformal field theory (Cardy [Car87], Di Francesco, Mathieu and Senechal [DiFMS96]).

Reviews of the Ising model include McCoy and Wu [MW73]. Cipra [Cip87] includes a highly accessible treatment of the 1-dimensional Ising model and the existence and nature of a phase transition for the 2-dimensional Ising model.  This is discussed in a more general setting in Stanley [Sta71].

\section{Extremal temperatures and the chromatic polynomial}

In addition to critical temperatures in $\left( {0,\infty } \right)$, and the straight-forward case of $\beta  = 0$ (i.e. infinite temperature), the extremal case of zero temperature, in both the ferromagnetic and antiferromagnetic models, is studied. Of particular interest, with respect to our theme of the interconnections between graph theory and statistical mechanics, is the equivalence of the zero temperature antiferromagnetic Potts model partition function and the chromatic polynomial.  We will first briefly mention the infinite temperature and the ferromagnetic cases, and then focus attention on the chromatic polynomial and its zeros.

In statistical physics, the free energy, $F = U - TS$, is minimized in a system in thermal equilibrium.  Note that at low temperatures, the minimization of the free energy becomes equivalent to the minimization of the internal energy $U$.  For positive $J$ this means that all spins take on the same value, while for negative $J$ it means that adjacent spins must have different values.  As $T$ increases, the minimization of the free energy increasingly means the maximization of the entropy, $S$.  As $T$ approaches infinity, the minimization of the internal energy plays a negligible role relative to the maximization of entropy.

At infinite temperature, i.e., when $K = 0$, we have $Z\left( {G;q,v} \right) = q^{\left| {V\left( G \right)} \right|} $, so $f = \ln q$. It is then possible to calculate the high-temperature Taylor series expansions for thermodynamic quantities as powers of $v = \exp (K) - 1$, so around $v = 0$.  These expansions are carried out around $v = 0$ because it is possible to systematically generate (e.g. via graphical techniques as in Nagle [Nag71] and Kim and Enting [KE79]), higher order terms in the Taylor series expansion about this point.

For the ferromagnetic Potts model at $T = 0$, the system exists in a completely ordered state, in which all spins have the same value, so that as $T \to 0$, and hence $K \to \infty $
, the partition function $P_1 \left( {G;q,\beta } \right) \to q\exp \left( {K\left| {E\left( G \right)} \right|} \right)$.  For a regular lattice whose vertices have uniform degree $\delta$, the reduced free energy per unit volume is then $\left| {V(G)} \right|^{-1}\ln q + \frac{\delta }{2}K$, so as the number of vertices goes to infinity, this gives the free energy per vertex as simply $F =  - (\delta /2){\rm{ }}J$. The low-temperature series expansion can then be expressed as a series in $\exp ( - K)$.  Reasonably accurate values for both the critical temperature and the critical exponents can then be extracted from the Taylor series expansion.

This now brings us to the zero-temperature antiferromagnetic case and the chromatic polynomial. We first recall that a proper coloring of a graph $G$ is an assignment of a color to each vertex of $G$ so that any two adjacent vertices receive different colors. The chromatic polynomial, $C(G;x)$, is a graph invariant that, when evaluated at a non-negative integer $x$, gives the number of ways to properly color the graph $G$ using $x$ colors.  Consider for a moment an edge $e$ of a graph $G$.  The number of ways to color $G-e$  (where there are no restrictions on the colors assigned to the endpoints of $e$) is equal to the number of ways to color $G$ (where the endpoints must have different colors) plus the number of ways to color $G/e$ (where the endpoints, now coalesced, must have the same color).  This means that the chromatic polynomial may be computed recursively as follows:
\[
 C(G - e;x) = C(G;x) + C(G/e;x) \textnormal{ or, } C(G;x) = C(G - e;x) - C(G/e;x),
\]
\[
\textnormal{ and } C(G;x) = x^n ,  \textnormal{ if } G \textnormal{ has } n \textnormal{ vertices and no edges.}
\]

Thus, the chromatic polynomial has a contraction/deletion reduction and hence by Theorem 3.1 must be an evaluation of the Tutte polynomial.  An argument almost identical to the proof of Theorem 3.2 shows that
\[
C(G;x) = ( - 1)^{\left| {V(G)} \right| - k(G)} x^{k(G)} t(G;1 - x,0).
\]

Some early reviews on the chromatic polynomial include those of Read [Rea68], Read and Tutte [RT88], and Biggs [Big93]. An extensive bibliography is available in Chia [Chi97], while Thomassen [Tho01] and Dong, Koh, and Teo [DKT05] give recent comprehensive treatments.

The connection between the Potts model and the chromatic polynomial occurs in the zero-temperature (limit as $T \rightarrow \infty$) antiferromagnetic  model using the $h_1 (\omega ) =  - J\sum\limits_{ij \in E(G)} {\delta (\sigma _i ,\sigma _j )} $ formulation of the Hamiltonian.  Since $J$ is negative in the antiferromagnetic model, minimal energy states are those that generate a maximum number of zeros in the summation, i.e. those in which every edge has endpoints with different spins.  If we think of the spins as colors, a minimum energy state is then just a proper coloring of the graph.

We give two different ways to understand the translation between the zero temperature Potts model and the chromatic polynomial.  One approach is to compare $C(G;x) = ( - 1)^{\left| {V(G)} \right| - k(G)} x^{k(G)} $ $t(G;1 - x,0)$  with the result of Theorem 3.2 that $P_1 \left( {G;q,\beta } \right) = q^{k(G)} v^{\left| {V(G)} \right| - k(G)} t\left( {G; (q + v)/v,v + 1} \right)$, where $v = \exp (J\beta ) - 1$.  We note that these are the same function precisely when $v=-1$, that is, when  $\beta = \infty$, which is exactly the zero-temperature model.  Another way to view this connection is by considering the summands of $P_1 (G;q,\beta ) = \sum {\exp \left( {\beta J\sum {\delta \left( {\sigma _i ,\sigma _j } \right)} } \right)} $.  As $T \rightarrow 0$, and hence  $\beta \rightarrow \infty$, a summand is 0 except precisely when $\sum{\delta(\sigma_i,\sigma_j)}=0$, in which case it is 1.  Thus $P_1(G)$ simply counts the number of proper colorings of $G$ with $q$ colors.

In the special case of the Potts antiferromagnetic model at $T = 0$, where the Potts model partition function reduces to the chromatic polynomial,  there are power series expansions of the various thermodynamic functions.  When $T = 0$, the model will be in one of its possible ground states.  Ground state entropy is a measure of the residual disorder in the system, and it can be nonzero for sufficiently large $q$.   In the infinite volume limit, the ground state entropy per vertex of the Potts antiferromagnetic model becomes

\[
S = \kappa \mathop {\lim }\limits_{n \to \infty } \frac{1}{{\left| {V\left( {G_n } \right)} \right|}}\ln \left( {C \left( {G_n ;q} \right)} \right).
\]

This is related to the ground state degeneracy per vertex, $W$, according to $S = \kappa \ln W$.  Two exact results are Lieb's calculation in [Lie67] giving $W = (4/3)^{3/2}$ for $q = 3$ on the square lattice (see also Baxter, Kelland and Wu [BKW76]) and Baxter's calculation of $W$ for general $q$ for the triangular lattice (see Baxter [Bax87], Wu [Wu82] and also Blote and Nightingale [BN82] and Baxter [Bax86]).  More generally, some calculations of Tutte polynomials for recursive families of graphs have been carried out in Shrock [Shr00], Chang and Shrock [CS00, 01a, 06], Chang, Salas and Shrock [CSS02], Change, Jacobsen, Salas and Shrock [CJSS04], and Jacobsen, Richard and Salas [JRS06].

A significant body of work has emerged in recent years devoted to clearing regions of the complex plane (in particular regions containing intervals of the real axis) of roots of the chromatic polynomial.  Results showing that certain intervals of the real axis and certain complex regions are free of zeros of chromatic polynomials include those given by Woodall [Woo92], Jackson [Jac93], Shrock and Tsai [ST97a,c], Thomassen [Tho97], Brown [Bro98], Sokal [Sok01b], Procacci, Scoppola, and Gerasimov [PSG03], Choe, Oxley, Sokal, and Wagner [COSW04], Borgs [Bor06], and  Fernandez and Procacci [FP].  One particular question concerns the maximum magnitude of a zero of a chromatic polynomial and of zeros comprising region boundaries in the $q$ plane as the number of vertices $|V| \rightarrow \infty $.  An upper bound is given in [Sok01b], depending on the maximal vertex degree.  There are families of graphs where both of these magnitudes are unbounded (see Read and Royle [RR91], Shrock and Tsai [ST97a, 98], Thomassen [Tho00], Bielak [Bie01], Brown, Hickman, Sokal and Wagner [BHSW01], and Sokal [Sok04]).  For recent discussions of some relevant research directions concerning zeros of chromatic polynomials and properties of their accumulation sets in the complex $q$ plane, as well as approximation methods, see, e.g., Shrock and Tsai [ST97b], Shrock [Shr01], Sokal [Sok01a, 01b], Chang and Shrock [CS01b], Chang, Jacobsen, Salas, and Shrock [CJSS04],  Choe, Oxley, Sokal, and Wagner [COSW04], Dong [Don04], Dong and Koh [DK04], and most recently Royle [Roya, b].

This study of the complex roots of chromatic polynomials extends previous work that traditionally focused on real roots, in particular, positive integer roots $q$ which correspond to a graph not being properly colorable with $q$ colors.

\section{Computational complexity connections}

The $q$-state Potts model partition function in Definition 2.2 involves a sum over all possible states of $G$.  If $G$ has $n$ vertices, then there are an exponential number, $q^n $, of states.  This immediately leads to the question of its computability. While realizing the Potts model partition function as an evaluation of the Tutte polynomial does not make it any easier to compute, it does enable the theory of one to inform the theory of the other and vice versa.  In fact, it was the computational complexity of the Tutte polynomial in general that showed rigorously that the Potts model partition function is likewise intractable.  The interplay of computational complexity results between the Tutte polynomial and the Potts model particularly illustrates the synergy between the two fields.

We first recall the basic notions of computational complexity.  A decision problem is one for which there is a yes or no answer, such as, can graph $G$ be colored using $k$ colors?  $P$ is the set of decision problems for which we can determine the answer in polynomial time in the size of the input, and $NP$ is the set of decision problems for which we can determine if a given answer is correct in polynomial time in the size of the input.  Whether or not $P=NP$ remains a famous open question, but there is a large class of $NP$-Hard problems for which finding a polynomial time algorithm for any one of them would automatically lead to polynomial time algorithms for all of them.  In practice, these $NP$-Hard problems are viewed to be intractable.  The set of $NP$-Complete problems are those decision problems in $NP$ that are known to be $NP$-Hard.  Analogously, the set of $\#P$-Complete problems is a complexity class consisting of counting problems (such as how many ways can a graph $G$ be colored using $k$ colors) that are similarly considered intractable.

Computational complexity results for the Potts model and the Tutte polynomial have built in alternation upon one another as the theory has evolved.  The 1990 paper of mathematicians Jaeger, Vertigan, and Welsh [JVW90] played a major role in this evolution.  In it, the authors note that the Ising model ($q=2$ Potts model) partition function can be reformulated as a tractable problem for planar graphs (referencing physicists Fisher [Fis66] and Kastelyn [Kas67]), but that it was shown by Jerrum [Jer87] to be $\#P$-Complete in general. (Also see Vertigan [Ver05].) The approach in [JVW90] focuses on the problem of evaluating the Tutte polynomial along hyperbolas of the form $(x-1)(y-1)=q$  where $q$ is any real number.   Note that if we let $x = (q + v)/v$  and $y=v+1$, and compare to Theorem 3.2, then these hyperbolas correspond precisely to the Potts model partition function for fixed positive integer values of $q$.  The conclusion of [JVW90] is that computing  the Tutte polynomial is $\#P$-Complete for general graphs, except when $q=1$  (which is trivial when viewed either in terms of the Potts model partition function or the Tutte polynomial), or when $q=2$  as discussed above, or for 9 special points, namely (1,1), (0,0), (-1,0), $( - i,i)$, $\left( {\exp (2\pi i/3),\exp (4\pi i/3)} \right)$ and their reflections about the line $y=x$, for which the Tutte polynomial has 'easy' enumerative interpretations.

Since the Tutte polynomial, and hence the Potts model partition function, is thus typically computationally intractable for arbitrary graphs and argument values, a natural question arises as to how well either might be approximated.  The answer is that in general approximating is provably difficult as well, but here again there is remarkable synergy between physics and mathematics, with results alternating between the fields.  We refer the reader to an excellent overview given by Welsh and Merino in Section VIII of [WM00], and to Alon, Frieze, and Welsh [AFW94, 95] for a more optimistic prognosis in the case of dense graphs.

There has also been an increasing body of work since the seminal results of Robertson and Seymour [RS83, 84, 86] addressing computational complexity questions for graphs with bounded tree-width (see Bodlaender's accessible introduction in [Bod93]).  A powerful aspect of this work is that many $NP$-Hard problems become tractable for graphs of bounded tree-width. Recent research includes a number of results both for the classical Tutte polynomial and also for the colored Tutte polynomial which encompasses the Potts model with variable interaction energies.  For example, Noble [Nob98] has shown that the Tutte polynomial may be computed in polynomial time (in fact requires only a linear number of multiplications and additions) for rational points on graphs with bounded tree width, and Makowsky [Mak05] and Traldi [Tra06] have extended this result to the colored Tutte polynomial.  Gimenez, Hlineny and Noy [GHN06] and Makowsky, Rotics, Averbouch and Godlin [MRAG06] provide similar results for bounded clique-width (a notion with significant computational complexity consequences analogous to those for bounded tree-width; see Oum and Seymour [OS06]).

While these computational complexity results for bounded tree- and clique-width are helpful in many instances, computing limits of the Potts model partition function as the number of vertices increases in an unbounded family of graphs remains an open question.  Even on families of lattices this may be problematic.  For example, Vertigan and Welsh [VW92] have shown that the Tutte polynomial is intractable away from $q=2$  even on planar bipartite graphs (except for certain trivial cases such as $q = 0$ or $1$), and Farr [Far06] shows that computing the number of colorings of induced subgraphs of even the square lattice is $\#P$-Complete.

Various approaches are used to circumvent this obstacle.  They include the Taylor series expansions previously discussed that provide powerful means of obtaining approximate information about the Potts model, as well as the computer simulations of the next section.   Additionally, calculation of the chromatic and Tutte polynomials may sometimes be achieved for a carefully chosen family of graphs where the iterative operation of the deletion-contraction property leads to a solvable closed set of linear equations.  Roughly speaking, the $\left( {m + 1} \right)^{\text{th}} $ member of such a family is constructed by gluing a particular subgraph to the $m^{\text{th}} $
 member (see Biggs, Damerill and Sands [BDS72]).  An example is a strip of a regular lattice of fixed width and variable length $m$.  The resulting Potts model partition function has the form of a finite sum of $m^{\text{th}}$ powers of a set of algebraic functions multiplied by certain coefficients.  These algebraic functions are the roots of a set of equations resulting from the iterative operation of the deletion-contraction theorem or equivalently, are the eigenvalues of a certain type of transfer matrix.  Some calculations of chromatic polynomials of recursive families of graphs include Biggs, Damerill and Sands [BDS72], Beraha and Kahane [BK79], Beraha, Kahane, and Weiss [BKW80], Read [Rea88], Salas [Sal90, 91] , Read and Royle [RR91], Shrock and Tsai [ST97a], Rocek, Shrock and Tsai [RST98], Shrock and Tsai [ST99], Shrock [Shr99, 01], Biggs and Shrock [BS99], Sokal [Sok00], Salas and Sokal [SS01], Chang and Shrock [CS01a, b], Biggs [Big01, 02a], Jacobsen, Salas and Sokal [JSS03], Jacobsen and Salas [JS01, 06].

\section{Monte Carlo simulations of the Potts model}

The computational intractability of the Potts model partition function has led to the development of Monte Carlo simulations for the model; see the texts of Newman and Barkema [NB99] and Landau and Binder [LB00] for additional background on the methods described below.  We illustrate the basic principles of this kind of simulation in the simplest, $q=2$, case of the Ising model on a square lattice, and then briefly mention some modifications leading to more sophisticated simulations.

Since complexes are often very large, with many different spin choices for their elements, the probability of a single state appearing out of the exponential number of states is nearly zero, but the macroscopic properties for many different states may be similar.  Therefore, the goal is to determine the average characteristics the system is likely to exhibit in the long run; i.e., we want to approximate the expected value of a macroscopic property when the system is in equilibrium.  In the case of the Ising model, we might want to determine the expected value of the magnetism at a given temperature.  The simulation must compute this expected value by averaging over a sufficiently large sample of states that correspond to an independent random sample of states from the Boltzmann distribution.  These states are generated through a Markovian random walk on the lattice.

It is preferable, but not necessary, to begin with an initial state that is characteristic of the temperature at which the properties of the system are being measured.  For example, if one were to start with an ordered spin configuration at a high temperature, then considerable computer time would be expended to "warm up" the simulation, while if one starts with a random spin configuration, much less time is spent reaching equilibration.  When the system is at equilibrium, the value of the macroscopic property of interest should stay within a fairly small range.  The simulation is generally run from a number of different initial configurations to ensure that the system has actually found the equilibrium value, rather than a locally stable value.  Since a state is clearly dependent on a few of the previous states in the random walk, an autocorrelation function is computed to determine the distance between samples taken in the random walk to ensure that the sample points are independent.  The necessary simulation length and corresponding statistical error can then be estimated in the typical manner for applications of the Central Limit Theorem.  As with all experiments, systematic error may occur and can be difficult to detect.

From the initial state, each vertex is visited in turn, and the program computes the probability ratio comparing the likelihood of the vertex changing its spin versus retaining its current spin.  This simulation captures the effect of temperature on the model, encoding the tendency of the system to move toward a lower energy state at low temperatures and remain agitated at high temperatures, as follows.  Recall that the probability of a state occurring is $\Pr \left( \varpi  \right) = (\exp \left( { - \beta h_i \left( \varpi  \right)} \right)/P_i \left( G \right)$, so the ratio of the probability of a new state $S_N $ to the probability of the current state $S_C $ is

\[
\frac{{\Pr (S_N )}}{{\Pr (S_C )}} = \frac{{\dfrac{{\exp \left( { - \beta h_i (S_N )} \right)}}{{\sum {\exp \left( { - \beta h_i (\omega )} \right)} }}}}{{\dfrac{{\exp \left( { - \beta h_i (S_C )} \right)}}{{\sum {\exp \left( { - \beta _i (\omega )} \right)} }}}} = \frac{{\exp \left( { - \beta h_i (S_N )} \right)}}{{\exp \left( { - \beta h_i (S_C )} \right)}} = \exp \left( {\frac{{h_i (S_C ) - h_i (S_N )}}{{\kappa T}}} \right).
\]

Note that this avoids computing the generally $NP$-Hard partition function, $P_i $.  Also note that since $S_N $
 differs from $S_C $ only in a change of spin at one vertex $v$, the computation of  $h_i \left( {S_N } \right)$ is exactly the same as that for $h_i \left( {S_C } \right)$ at every edge except for those incident with $v$.

In the commonly used Metropolis Monte Carlo algorithm, if the new state has lower energy than the current state, $h_i \left( {S_N } \right) \le h_i \left( {S_C } \right)$, the algorithm changes the system from state $S_C $  to state $S_N $.  However, if $h_i \left( {S_C } \right) < h_i \left( {S_N } \right)$, the program compares
 \[
\frac{{\Pr (S_N )}}{{\Pr (S_C )}} = \exp \left( {\frac{{h_i (S_C ) - h_i (S_N )}}{{\kappa T}}} \right)
\]
 to a random number $r$ with $0 \le r \le 1$ and changes state if $r < \Pr (S_N )/\Pr (S_C )$.  At high temperatures, this ratio will be nearly 1 regardless of the Hamiltonians, so spins will continue changing with negligible preference for lower energy states.  On the other hand, if the temperature is quite low, the system strongly favors low energy states.

Although the behavior of the model is clear at very high or very low temperatures, it is less apparent what happens at midrange temperatures.  One of the fundamental questions for the Potts model on a regular $d$-dimensional lattice is determining the critical temperature $T_C $ for a phase transition.  In Monte Carlo simulations, temperatures near the critical value can cause computational challenges, due to the increased statistical error in that range as well as the increase in the autocorrelation time. Nevertheless, simulations are an important tool in the study of the Potts model, as few exact analytic results are known.

There are various refinements of this basic model leading to more sophisticated simulations.  For example, there is no need for the underlying graph to be a square lattice. It can be a different regular lattice, such as triangular or honeycomb in two dimensions; a cubic, face-centered cubic, or body-centered cubic in three dimensions. More generally, it can be any graph appropriate to the application, even a complete graph if every site interacts with every other, although dense  and/or irregular graphs can present programming challenges.  The simulation can be extended to larger numbers of spins by computing the relative probabilities comparing the current spin at a vertex with each of the other possible spins.  In the heat-bath algorithm, the probability ratios are normalized so they sum to 1, and then each is assigned a proportional segment of the unit interval.  A random number is generated in the unit interval, and the spin is changed according to the segment that contains the random number.
There are also useful techniques for improving the speed of these simulations, including clustering methods.  Clusters of locally aligned states slow down the simulation, since the likelihood of a "flip" occurring is very low.  Therefore, the simulation will spend extended periods in the same state.  The Wolff algorithm improves the running time by flipping these clusters of like spins together instead of considering them one by one. See Wang, Kozen, and Swendsen [WKS02], Deng, Garoni, and Sokal [DGS07], and Deng, Garoni, Machta, Ossola, Polin, and Sokal [D+07].

Other significant modifications and variations of the model include inclusion of an external magnetic field, next-nearest-neighbor interactions, edge-dependent interaction energies, and additional terms in the Hamiltonian.  In realistic physics studies, one also must often consider the effects of disorder in actual materials, such as vacancies and impurities in a crystal lattice.  We give further discussion of some of these variations in the next section.

\section{Why is this model attracting so much attention?}

Besides its intrinsic mathematical interest, the Potts model, in many variations, is increasingly applicable to a wide variety of complex systems where local interactions can predict global behavior.  This is particularly true as computing power has enabled increasingly powerful and predictive simulations and as researchers have found sophisticated modifications of the model to more closely mimic the behaviors of various systems. The popularity of the Potts model is roughly indicated by a recent Google Scholar search for ``Potts model'' producing over 63,000 hits.  We give a sample of applications here, just to demonstrate the scope of this theory.  In these examples, the Hamiltonian is extended to encode forces in addition to simple nearest-neighbor interactions, but the probability distribution, and hence partition function, is still generally defined analogously to that of Definition 2.2.

The original magnetism application addressed by Ising considers the overall ferromagnetic (``normal'' magnetism) behavior of a lattice where the two possible spins at each position are positive and negative.  The energy of the system is minimized if all points on the lattice have the same spin, while, to maximize entropy, all states should be equally likely (which would strongly favor nonmagnetic states).  The Boltzmann distribution quantifies the relative importance of energy and entropy in determining the likelihood of a given state in terms of temperature (assuming the system is at the same temperature as its surrounding environment).  The standard Hamiltonian is given in Definition 2.1, but is sometimes extended to include an external magnetic field:
\[
h(\omega ) =  - J\sum\limits_{ ij  \in E(G)} {\delta (\sigma _i ,\sigma _j )}  - H\sum\limits_i {\sigma _i }.
\]

Simulations and series expansions are used in higher dimensions to determine the phase transition temperature below which the system exhibits a nonzero spontaneous magnetization and above which this magnetization vanishes.   See Stanley [Sta71], McCoy and Wu [MW73], Chandler [Cha87], Plischke and Bergesen [PB06] for the Ising model, and Stanley [Sta71], Carter [Car01], Plischke and Bergesen [PB06] for further information on the Boltzmann distribution.

Sanyal and Glazier [SG06] employ the Potts model to simulate foam flow and investigate instabilities, i.e., velocities at which larger bubbles start to flow faster than smaller bubbles.  Two adjacent lattice points have the same spin if and only if they are part of the same bubble, hence the number of spins $Q$ is extremely large.  The Hamiltonian is:
\[
H = \sum\limits_{ij \in E\left( G \right)} {J\left( {1 - \delta \left( {\sigma _i ,\sigma _j } \right)} \right)}  + \lambda \sum\limits_{n = 1}^Q {\left( {a_n  - A_n } \right)^2 },
\]
where $J$ is the coupling strength at the boundary between two bubbles, an is the actual area of the bubble, $A_n$ is the area the bubble would have if it were not subjected to external forces, and $\lambda$ is the strength of the area constraint on the bubble (based on the compressibility of the gas).  At each step in the simulation, the spin at a randomly selected lattice site is considered; if that site is along the boundary with another bubble, a switch in spin of that lattice site to the neighboring bubble is considered (and accepted with some probability).

Turner and Sherratt [TS02] use an extension of the Potts model to study cellular malignancy growth.  They are particularly interested in the impact of the relative strength of a few factors known to affect cell behavior.  The Hamiltonian for their model is:
\[
H = \sum\limits_{ij} {\sum\limits_{i'j'} {J_{\sigma _{ij} ,\sigma _{i'j'} } \{ 1 - \delta _{\tau (\sigma _{ij} ),\tau (\sigma _{i'j'} )} \}  + \sum\limits_\sigma  {\lambda (v_\sigma   - V_T )^2 } } } .
\]

Here $ij$ indexes the lattice point in the $i^{th}$ row and $j^{th}$ column of a  two dimensional grid.  The spin value $\sigma_{ij}$  records which of the many cells contains the $ij^{th}$ lattice point.  The first term of the Hamiltonian encodes at the interaction energy between a cell and its eight nearest neighbors on the lattice.  Adjacent lattice sites with the same spin value  $\sigma_{ij}$  represent a single cell, so there will be no interaction in that case.  Otherwise the interaction energy depends on the cell type $\tau(\sigma_{ij})$, which may be normal or malignant.  The second term models the energy required for a cell $\sigma $ to maintain a volume different from its natural volume in the absence of external forces, similar to the foam example above.  Instead of temperature, the $\beta $  in this application corresponds to a diffusion coefficient affecting the random motility of the cells.  Turner and Sherratt [TS02] further extend the Potts model by allowing cells to replicate, hence changing the lattice, during the simulation.  In doing so, they take into account the interaction of a cell with the extracellular (protein) matrix.

Nobel laureate Thomas Schelling published a seminal paper titled \emph{Dynamic models of segregation} in 1971 [Sch71] that considers the possibility of micro-motive explanations for racial segregation (in addition to organized and economic explanations).  The premise is that individual decisions to avoid minority status (or to require being in a minority of some minimum size) could lead to the macro-effect of segregation.  Schelling places vacancies, stars, and zeros randomly on a checkerboard and then iteratively considers the 'happiness' of the stars and zeros with their local neighborhoods, moving an unhappy star or zero to the nearest vacant spot that meets their happiness criteria.  Meyers-Ortmanns [Mey03] models a similar premise to Schelling's (that micro-motive explanations can lead to immigrant ghettos) with a more Potts-like model where the Hamiltonian measures the happiness of individuals with their neighbors, the temperature is viewed as a social temperature where warmer temperatures reflect facilitation of integration and assimilation, and at each step in the simulation two neighbors are able to exchange places with a probability based on the likelihood of the new state with respect to the current state.  Schulze [Sch05] extends Meyers-Ortmanns to address up to seven different ethnic groups.

\section{Some active areas at the interface of combinatorics and statistical mechanics}

While there is intense interest in the Potts model, not just from the physics community, but throughout the sciences, its properties are of intrinsic combinatorial interest as well, and clearly much work remains to be done in exploring, developing, and extending them.  For readers who would like to pursue broad perspectives and further background on the topics mentioned here, overviews of the Tutte polynomial may be found in Brylawski [Bry80], Tutte [Tut84],  Brylawski and Oxley [BO92], Biggs [Big93], Welsh [Wel93, 99], Bollob\'{a}s [Bol98], and Farr [Far07], with the latter three also discussing the Potts model from a mathematical point of view.  Some works containing reviews of the Potts model from physics and mathematical physics viewpoints include those of Baxter [Bax82], Wu [Wu82, 88], Cipra [Cip87], Martin [Mar91], Welsh and Merino [WM00], Shrock [Shr01], Chang, Jacobsen, Salas, and Shrock [CJSS04], Sokal [Sok00, 05], and Farr [Far07a].

In addition to the areas already discussed, researchers are also currently very interested in a number of related areas.  New computational techniques for relevant combinatorial polynomials, particularly those that can be applied to regular structures such as lattices, are always sought. Of independent, yet mutual, interest for both graph theory and statistical mechanics is the determination and physical interpretation of the zeros of Tutte polynomial and its cognates.  The effects of boundary conditions on lattices, including periodicity leading to toroidal and other topologies, are an important area of study, as are the connections to knot theory.  Improved Markov-chain Monte Carlo methods are in constant demand, as are further computational complexity results.  Recent work also includes the use of methods from statistical mechanics in combinatorial enumeration. Space prohibits providing an exhaustive list, but a few examples not previously mentioned in these areas include Chang and Shrock [CS01c]; Sokal [Sok01a]; Wu and Wang [WW01]; Woodall [Woo02]; Traldi [Tra02, 06]; Bonin and de Meir [BdMN03];  Jackson [Jac03]; Read [Rea03]; Biggs [Big04], Royle and Sokal [RS04]; Morris and Peres [MP05]; Makowsky [Mak05]; Wu [Wu05, 06]; Makowsky, Rotics, Averbouch, and Godlin [MRAG06]; Oum and Seymour [OS06]; Farr [Far07b]; Martinelli, Sinclair, and Weitz [MSW07]; Jerrum [Jer07], and Jacobsen and Salas [JS07].

We close with an observation that, although far from original, hopefully may encourage combinatorialists to reach across the dark spaces between disciplines and engage in research related to statistical mechanics.  Theoretical physicists have produced a wealth of information about phase transitions and critical phenomena leading to well-supported assertions, many of which still need rigorous mathematical treatment or lead to questions of intrinsic mathematical interest. The result is a ready supply of appealing mathematical problems.  This is especially true for combinatorialists, into whose domain many of these problems from statistical mechanics naturally fall.

\textbf{\textit{Acknowledgements}}

We thank Alain Brizard, William Karstens, Lorenzo Traldi, Peter Winkler, Thomas Zaslavsky, and especially Alan Sokal for a number of informative conversations.

J. Ellis-Monaghan and R. Shrock thank the Isaac Newton Institute for Mathematical Sciences of Cambridge University for hospitality during the time when part of this work was completed.  The Combinatorics and Statistical Mechanics Programme of the Newton Institute contains a number of valuable reviews and research presentations on the subjects discussed here, and these are available online at [New08].

Support was provided by the National Security Agency and by the Vermont Genetics Network through Grant Number P20 RR16462 from the INBRE Program of the National Center for Research Resources (NCRR), a component of the National Institutes of Health (NIH).  This paper's contents are solely the responsibility of the authors and do not necessarily represent the official views of NCRR or NIH.

\section{Bibliography}

[AFW94] Alon, N.; Frieze, A.; Welsh, D.  Polynomial time randomized approximation schemes for the Tutte polynomial of dense graphs.  35th Annual Symposium on Foundations of Computer Science (Santa Fe, NM, 1994), IEEE Comput. Soc. Press, Los Alamitos, CA, (1994), 24-35.

[AFW95] Alon, N.; Frieze, A.; Welsh, D.  Polynomial time randomized approximation schemes for Tutte-Gröthendieck invariants:  the dense case.  Random Stuctures Algorithms 6 (1995), no. 4, 459-478.

[Ale75] Alexander, S. Lattice gas transition of He on Grafoil: a continuous transition with cubic terms. Phys. Lett. A 54 (1975), 353-354.

[AT43]  Ashkin, J.; Teller, E. Statistics of two-dimensional lattices with four components.  Phys. Rev. 64 (1943), no. 5 and 6.

[Bax73] Baxter, R.J. Potts model at the critical temperature. J. Phys. C 6 (1973), L445-L448.

[Bax82] Baxter, R.J. Exactly Solved Models in Statistical Mechanics. Academic Press, New York, 1982.

[Bax82a] Baxter, R.J. Critical square-lattice Potts model. Proc. Roy. Soc. London A 383 (1982), 43-54.

[Bax86] Baxter, R.J.  $q$ Colourings of the triangular lattice. J. Phys. A 19 (1986), 2821-2839.

[Bax87] Baxter, R.J.  Chromatic polynomials of large triangular lattices. J. Phys. A 20 (1987),  5241-5261.

[BdMN03]    Bonin, J.; de Mier, A.; Noy, M. Lattice path matroids: enumerative aspects and Tutte polynomials. J. Combin. Theory Ser. A 104 (2003), no. 1, 63-94.

[BDS72] Biggs, N.L.; Damerill, R.M.; Sands, D.A.  Recursive families of graphs.  J. Combin. Theory 12 (1972), 123-131.

[BHSW01]    Brown, J.I.; Hickman, C.A.; Sokal, A.D.; Wagner, D.G.  On the chromatic roots of generalized theta graphs. J. Combin. Theory B 83 (2001), 272-297.

[BK79]  Beraha, S.; Kahane, J.  Is the four-color theorem almost false? J. Combin. Theory 27 (1979), 1-12.

[Bie01] Bielak, H. Roots of chromatic polynomials. Discrete Math. 231 (2001), 97-102.

[Big93] Biggs, N.L. Algebraic Graph Theory (2nd ed.). Cambridge University Press, Cambridge (1993).

[Big01] Biggs, N.L. Matrix method for chromatic polynomials. J. Combin. Theory B 82 (2001), 19-29.

[Big02a]    Biggs, N.L. Chromatic polynomials for twisted bracelets. Bull. London Math. Soc. 34 (2002), 129-139.

[Big02b]    Biggs, N.L. Equimodular curves. Discrete Math. 259 (2002), 37-57.

[Big04] Biggs, N.L. Specht modules and chromatic polynomials. J. Combin. Theory, Ser. B 92 (2004), 359-377.

[BK04]  Biggs, N.L.; Klin, M.H.; Reinfeld, P. Algebraic methods for chromatic polynomials. Eur. J. Combin. 25 (2004), 147-160.

[BKW76] Baxter, R.; Kelland, S.B.; Wu, F.Y.  Equivalence of the Potts or Whitney polynomial to an ice model. J. Phys. A 9 (1976), 397-406.

[BKW80] Beraha, S.; Kahane, J.; Weiss, N. Limits of chromatic zeros for some families of maps.  J. Combin. Theory, Ser. B 28 (1980), 52-65.

[BN82]  Blote, H.W.J.; Nightingale, M.P. Critical behaviour of the two-dimensional Potts model with a continuous number of states: a finite-size scaling analysis. Physica A 112 (1982), 405-465.

[BO92]  Brylawski, T.; Oxley, J. The Tutte polynomial and its applications. Matroid Applications, Encyclopedia Math. Appl. 40, Cambridge Univ. Press, Cambridge. 1992, pp. 123-225.

[Bod93] Bodlaender, H.L.  A tourist guide through treewidth.  Acta Cybernet 11 (1993), 1-21.

[Bol98] Bollob\'{a}s, B.  Modern Graph Theory.  Graduate Texts in Mathematics.  Springer-Verlag New York, Inc.,  New York, New York, 1998.

[Bor06] Borgs, C. Absence of zeros for the chromatic polynomial on bounded degree graphs. Combin. Probab. Comput. 15 (2006), no. 1-2, 63-74.

[BR99]  Bollob\'{a}s, B.; Riordan, O.   A Tutte polynomial for coloured graphs. Recent Trends in Combinatorics (Mátraháza, 1995). Combin. Prob. Comp. 8 (1999), 45-93.

[Bro98] Brown, J.  On the roots of chromatic polynomials. J. Combin. Theory 72 (1998), 251-256.

[Bry80]     Brylawski, T. The Tutte polynomial. Proceedings of the Third International Mathematical Summer Centre (1980), 125-275.

[BS99]  Biggs, N.L.; Shrock, R. $T=0$ partition functions for Potts antiferromagnets on square lattice strips with (twisted) periodic boundary conditions. J. Phys. A (Letts.) 32 (1999), L489-L493.

[BW99]  Brightwell, G.R.; Winkler, P. Graph homomorphisms and phase transitions. J. Combin. Theory Ser. B 77 (1999), no. 2, 221-262.

[Car87] Cardy, J. Conformal invariance. Phase Transitions and Critical Phenomena 11, eds. C. Domb and J. Lebowitz, Academic Press, New York 1987, pp. 55-126.

[Car01] Carter, A.H. Classical and Statistical Thermodynamics. Prentice Hall, Inc., Upper Saddle River, NJ, (2001) 238-241.

[Cha87] Chandler, D. Introduction to Modern Statistical Mechanics.  Oxford University Press,  New York, New York, (1987).

[Chi97] Chia, G.L. A Bibliography on Chromatic Polynomials. Chromatic Polynomials and Related Topics (Shanghai, 1994). Discrete Math. 172 (1997), no. 1-3, 175-191.

[CHW96] Chen, C.-N.; Hu, C.-K.; Wu, F.-Y. Partition function zeros of the square lattice Potts model. Phys. Rev. Letts. 76 (1996), 169-172.

[Cip87] Cipra, B.A.  An introduction to the Ising model.  Amer. Math. Monthly  94 (1987), no. 10, 937-959.

[CJSS04]    Chang, S.-C.; Jacobsen, J.; Salas, J.; Shrock, R. Exact Potts model partition functions for strips of the triangular lattice. J. Stat. Phys. 114 (2004), 768-823.

[COSW04]    Choe, Y.B.; Oxley, J.G.; Sokal, A.D.; Wagner, D.G. Homogeneous multivariate polynomials with the half-plane property. Adv. in Appl. Math. 32, (2004), no. 1-2, 88-187.

[CS00]  Chang, S.-C.; Shrock, R. Exact Potts model partition functions on strips of the triangular lattice. Physica A 286 (2000), 189-238.

[CS01a] Chang, S.-C.; Shrock, R. Ground state entropy of the Potts antiferromagnet on strips of the square lattice. Physica A 290 (2001), 402-430.

[CS01b] Chang, S.-C.; Shrock, R.  Exact Potts model partition functions on wider arbitrary-length strips of the square lattice. Physica A 296 (2001), 234-288.

[CS01c] Chang, S.-C.; Shrock, R. Zeros of Jones polynomials for families of knots and links, Physica A 301 (2001), 196-218.

[CS06]  Chang, S.-C.; Shrock, R.  Partition function zeros of a restricted Potts model on lattice strips and effects of boundary conditions. J. Phys. A 39 (2006), 10277-10295.

[CSS02] Chang, S.-C.; Salas, J.; Shrock, R. Exact Potts model partition functions for strips of the square lattice. J. Stat. Phys. 107 (2002), 1207-1253.

[D+07]  Deng, Y.; Garoni, T.; Machta, J.; Ossola, G.; Polin, M., and Sokal, A.D. Dynamic critical behavior of the Chayes-Machta-Swendsen-Wang algorithm. Preprint, arXiv:0705.2751.

[DGS07] Deng, Y.; Garoni, T.; Sokal, A.D. Critical speeding-up in the local dynamics of the random cluster model.  Phys. Rev. Letts. 98 (2007), 230602-1 to 23602-4.

[DiFMS96]   Di Francesco, P.;  Mathieu, P.; and Senechal, D. Conformal field theory. Springer, New York, (1996).

[DK04]  Dong, F.M.; Koh, K.M., On upper bounds for real roots of chromatic polynomials. Discrete Math. 282 (2004), 95-101.

[DKT05] Dong, F.M.; Koh, K.M; Teo, K.L. Chromatic polynomials and chromaticity of graphs. World Scientific Publishing Co. Pte. Ltd., Hackensack, NJ, 2005.

[Don04] Dong, F.M. The largest non-integer zero of chromatic polynomials of graphs with fixed order. Discrete Math. 282 (2004), 103-112.

[DS79]  Domany, E.; Schick, M. Classification of continuous order-disorder transitions in adsorbed monolayers II. Phys. Rev. B 20 (1979), 3828-3836.

[E-MT06]    Ellis-Monaghan, J.; Traldi, L.  Parametrized Tutte polynomials of graphs and matroids. Combin., Prob. Comp. 15 (2006), 835-854.

[Far06] Farr, G.  The complexity of counting colourings of subgraphs of the grid.  Comb. Probab. Comp. 15 (2006) 377-383.

[Far07a]    Farr, G. Tutte-Whitney polynomials: some history and generalizations. Combinatorics, Complexity, and Chance: a Tribute to Dominic Welsh, Oxford University Press, Oxford, 2007.

[Far07b]    Farr, G.E. On the Ashkin-Teller model and Tutte-Whitney functions. Combin. Probab. Comput. 16 (2007), no. 2, 251-260.

[Fis65] Fisher, M.E. The nature of critical points. Lectures in Theoretical Physics, vol. VIIC, University of Colorado Press, Boulder, CO, 1965, pp. 1-159.

[Fis66] Fisher, M.E.  On the dimer solution of planar Ising models.  J. Math. Phys. 7 (1966), 1776-1781.

[Fis74] Fisher, M.E. The renormalization group in the theory of critical phenomena. Rev. Mod. Phys. 46 (1974), 597-616.

[FK72]  Fortuin, C.M.; Kasteleyn, P.W. On the random cluster model. Physica (Amsterdam) 57 (1972), 536-564.

[FP]    Fernandez, R.; Procacci, A. Regions without complex zeros for chromatic polynomials on graphs with bounded degree. Preprint.

[GHN06] Gimenez, O.; Hlineny, P.; Noy, M. Computing the Tutte polynomial on graphs of bounded clique-width. SIAM J. Discrete Math. 20 (2006), 932-946.

[Isi25] Ising, E. Beitrag zur Theorie des Ferromagnetismus. Zeitschrift für Physik  31 (1925),  253-258.

[Jac93] Jackson, B. A zero-free interval for chromatic polynomials of graphs. Combin. Probab. Comput. 2 (1993), no. 3, 325-336.

[Jac03] Jackson, B. Zeros of chromatic and flow polynomials of graphs. J. Geom. 76 (2003), no. 1-2, 95-109.

[Jer87]     Jerrum, M.R.  2-dimensional monomer-dimer systems are computationally intractable.  J. Statist. Phys. 48 (1987), 121-134.

[Jer07] Jerrum, M.R. Approximating the Tutte polynomial. Combinatorics, complexity, and chance, Oxford Lecture Ser. Math. Appl. 34, Oxford Univ. Press, Oxford, (2007), 144-161.

[JRS06] Jacobsen, J.L.; Richard, J.-F.; Salas, J. Complex-temperature phase diagram of Potts and RSOS models.  Nucl. Phys. B 743 (2006), 153-206.

[JS07]  Jacobsen, J.L.; Salas, J. Phase diagram of the chromatic polynomial on a torus. Nucl. Phys. B 783, (2007), 238-296.

[JS01]  Jacobsen, J.L.; Salas, J.  Transfer matrices and partition function zeros for antiferromagnetic Potts models. II. Extended results for square lattice chromatic polynomial. J. Stat. Phys. 104 (2001), 701-723.

[JS06]  Jacobsen, J.L.; Salas, J.  Transfer matrices and partition function zeros for antiferromagnetic Potts models. IV. Chromatic polynomials with cyclic boundary conditions.  J. Stat. Phys. 122 (2006), 705-760.

[JSS03] Jacobsen, J.L.; Salas, J.; Sokal, A.D.  Transfer matrices and partition function zeros for the antiferromagnetic Potts model. III. Triangular lattice chromatic polynomial. J. Stat. Phys. 112 (2003), 921-1017.

[JVW90] Jaeger, F.; Vertigan, D.L.; Welsh, D.J.A. On the computational complexity of the Jones and Tutte polynomials.  Math. Proc. Cambridge Philos. Soc. 108 (1990), no. 1, 35-53.

[Kas67] Kasteleyn, P.W.  Graph theory and crystal physics.  Graph Theory and Theoretical Physics, ed. F. Harary,  Academic Press, 1967, pp. 43-110.

[KC01]  Kim, S.-Y.; Creswick, R. Density of states, Potts zeros and Fisher zeros of the $q$-state Potts model for continuous $Q$.  Phys. Rev. E63 (2001), 066107.

[KE79]  Kim, D.; Enting, I.G. The limit of chromatic polynomials. J. Combin Theory, Ser. B 26 (1979), no. 3, 327-336.

[LB00]  Landau, D.P.; Binder, K. A Guide to Monte Carlo Simulations in Statistical Physics. Cambridge University Press, Cambridge, 2000.

[Lie67] Lieb, E.H.  Residual entropy of square ice. Phys. Rev. 162 (1967) 162-172.

[LY52]  Lee, T.-D.; Yang, C.-N. Statistical theory of equations of state and phase transitions. II. Lattice gas and Ising model.  Phys. Rev. (2) 87 (1952), 410-419.

[Mak05] Makowsky, J.A. Colored Tutte polynomials and Kauffman brackets for graphs of bounded tree width. Discrete Appl. Math. 145 (2005), 276-290.

[Mar91] Martin, P.  Potts Models and Related Problems in Statistical Mechanics. World Scientific, Singapore, 1991.

[Mey03] Meyer-Ortmanns, H. Immigration, integration and ghetto formation. Internat. J. Modern Phys. C 14 (2003), no. 3, 311-320.

[MM86]  Martin, P.; Maillard, J.-.M. Partition function zeros for the triangular three-state Potts model. J. Phys. A 19 (1986),  L547-L550.

[MP05]  Morris B.; Peres, Y. Evolving sets, mixing and heat kernel bounds. Probab. Theory Related Fields 133 (2005), no. 2, 245-266.

[MRAG06]    Makowsky, J.A.; Rotics, U.; Averbouch, I.; Godlin, B. Computing graph polynomials on graphs of bounded clique-width. WG 2006, Lecture Notes in Computer Science 4271. Springer-Verlag New York Inc.  New York, New York. 2006, pp. 191-204.

[MS96]  Matveev, V.; Shrock, R. Complex-temperature singularities in Potts models on the square lattice. Phys. Rev. E54 (1996), 6174-6185.

[MSW07] Martinelli, F.; Sinclair, A.; Weitz, D. Fast mixing for independent sets, colorings and other models on trees. Random Structures Algorithms, 31 (2007), no. 2, 134-172.

[MW73]  McCoy, B.; Wu, T.T. The Two-Dimensional Ising Model. Harvard University Press, Cambridge, 1973.

[Nag71] Nagle, J.F. A new subgraph expansion for obtaining coloring polynomials for graphs. J. Combin. Theory, Ser. B 10 (1971), 42-59.

[NB99]  Newman, M.E.J.; Barkema, G.T. Monte Carlo Methods in Statistical Physics. Clarendon Press, Oxford, (1999).

[New08]  Combinatorics and Statistical Mechanics Programme, January-June 2008, Isaac Newton Institute for Mathematical Sciences, Cambridge University.  Videos of presentations and other resources available at http://www.newton.cam.ac.uk/webseminars/pg+ws/2008/csm/.

[Nob98]     Noble, S.D. Evaluating the Tutte polynomial for graphs of bounded tree-width. Combin. Prob. Comp.  7 (1998), no. 3, 307-321.

[Ons44] Onsager, L. Crystal statistics. I. A two-dimensional model with an order-disorder transition. Phys. Rev. (2) 65 (1944), 117-149.

[OS06]  Oum, S.; Seymour, P. Approximating clique-width and branch-width. J. Combin. Theory Ser. B 96 (2006), no. 4, 514–528.

[PB06]  Plischke, M.; Bergesen, B. Equilibrium Statistical Mechanics, 3rd ed. World Scientific, Singapore, 2006.

[Pot52] Potts, R.B. Some generalized order-disorder transformations. Proc. Cambridge Philos. Soc. 48 (1952), 106-109.

[PSG03]     Procacci, A.; Scoppola, B.; Gerasimov, V. Potts model on infinite graphs and the limit of chromatic polynomials. Communications in Mathematical Physics 235 (2003), no. 2, 215-231.

[Rea68] Read, R.C.  Introduction to chromatic polynomials.  J. Combin. Theory  4 (1968), 52-71.

[Rea88] Read, R.C.  A large family of chromatic polynomials. Proceedings of the Third Caribbean Conference on Combinatorics and Computing (1981), 23-41.

[Rea03] Read, R.C. Chain polynomials of graphs, Discrete Math. 265 (2003), 213-235.

[Roya]  Royle, G. Planar triangulations with real chromatic roots arbitrarily close to four. Preprint, math.CO.0511205.

[Royb]  Royle, G. Graphs with chromatic roots in the interval (1,2). Preprint arXiv:0704.2264.

[RR91]  Read, R.C.; Royle, G. Chromatic roots of families of graphs. Graph Theory, Combinatorics, and Applications vol. 2, eds. Y. Alavi et al., Wiley, New York, (1991), 1009-1029.

[RS83]  Robertson, N.; Seymour, P.D.  Graph minors. I. Excluding a forest. J. Combin. Theory Ser. B 35 (1983), 39-61.

[RS84]  Robertson, N.; Seymour, P.D.  Graph minors. III. Planar tree-width. J. Combin.  Theory Ser. B 36 (1984), 49-64.

[RS86]  Robertson, N.; Seymour, P.D.  Graph minors. II. Algorithmic aspects of tree-width.  J Algorithms 7 (1986), 309-322.

[RS04]  Royle, G.; Sokal, A. The Brown-Colbourn conjecture on zeros of reliability polynomials is false. J. Combin. Theory Ser. B 91 (2004), no. 2, 345-360.

[RST98] Rocek, M.; Shrock, R.; Tsai, S.-H. Chromatic polynomials for strip graphs and their asymptotic limits. Physica A 252 (1998), 505-546.

[RT88]  Read, R.C.; Tutte, W.T. Chromatic polynomials. Selected Topics in Graph Theory, vol. 3. eds. L. W. Beineke and R. J. Wilson, Academic, New York, (1988).

[Sal90] Saleur, H. Zeros of chromatic polynomials: a new approach to Beraha conjecture using quantum groups. Commun. Math. Phys. 132 (1990), 657-679.

[Sal91] Saleur, H. The antiferromagnetic Potts model in two dimensions: Berker-Kadanoff phase, antiferromagnetic transitions, and the role of Beraha numbers. Nucl. Phys. B 360 (1991), 219-263.

[Sch71] Schelling, T.C. Dynamic models of segregation.  Journal of Mathematical Sociology 1 (1971), 143-186.

[Sch05] Schulze, C. Potts-Like model for ghetto formation in multi-cultural societies. Internat. J. Modern Phys. C 16 (2005), Issue 03, 35-355.

[SG06]  Sanyal, S.;  Glazier, J.A. Viscous instabilities in flowing foams:  A cellular Potts model approach.  Journal of Statistical Mechanics. doi:10.1088/1742-5468/2006/10/P10008, (2006).

[Shr99] Shrock, R. $T=0$ Partition functions for Potts antiferromagnets on Mobius strips and effects of graph topology. Phys. Letts. A 261 (1999), 57-62.

[Shr00] Shrock, R. Exact Potts model partition functions for ladder graphs.  Physica A283(2000), 388-446.

[Shr01] Shrock, R. Chromatic polynomials and their zeros and asymptotic limits for families of graphs. Discrete Math. 231 (2001), 421-446.

[Sok00] Sokal, A.D.  Chromatic polynomials, Potts models and all that.  Physica A279 (2000), 324-332.

[Sok01a]    Sokal, A.D. A personal list of unsolved problems concerning lattice gases and antiferromagnetic Potts models.  Markov Processes and Related Fields 7 (2001), 21-38.

[Sok01b]    Sokal, A.D.  Bounds on the complex zeros of (di)chromatic polynomials and Potts-model partition functions.  Combin. Probab. Comput. 10 (2001), no. 1, 41-77.

[Sok04] Sokal, A.D. Chromatic roots are dense in the whole complex plane. Combin. Probab. Comput. 13 (2004), 221-261.

[Sok05] Sokal, A.D. The multivariate Tutte polynomial (alias Potts model) for graphs and matroids.  Surveys in Combinatorics, edited by Bridget S. Webb.  Cambridge University Press, 2005, pp. 173-226.

[SS97]  Salas, J.; Sokal, A.D. Absence of phase transitions for the antiferromagnetic Potts model via the Dobrushin uniqueness theorem.  J. Stat. Phys. 86 (1997), 551-579.

[SS01]  Salas, J.; Sokal, A.D. Transfer matrices and partition-function zeros for antiferromagnetic Potts models: I. General theory and square lattice chromatic polynomial. J. Stat. Phys. 104 (2001), 609-699.

[ST97a] Shrock, R.; Tsai, S.-H. Asymptotic limits and zeros of chromatic polynomials and ground state entropy of Potts antiferromagnets. Phys. Rev. E 55 (1997), 5165-5179.

[ST97b] Shrock, R.; Tsai, S.-H. Ground state entropy of Potts antiferromagnets: bounds, series, and Monte Carlo measurements.  Phys. Rev. E 56 (1997), 2733-2737.

[ST97c] Shrock, R.; Tsai, S.-H. Families of graphs with Wr(G,q) Functions That Are Nonanalytic at 1/q=0. Phys. Rev. E 56 (1997), 3935-3943.

[ST98]  Shrock, R.; Tsai, S.-H. Ground state entropy of  Potts antiferromagnets: cases with noncompact W boundaries having multiple points at $1/q=0$. Physica A 259 (1998), 315-348.

[ST99]  Shrock, R.; Tsai, S.-H.  Ground state degeneracy of Potts antiferromagnets on 2D lattices: approach using infinite cyclic strip graphs.  Phys. Rev. E 60 (1999), 3512-3515.

[Sta71]     Stanley, H.E. Introduction to Phase Transitions and Critical Phenomena. Oxford University Press, Oxford (1971).

[Tho97] Thomassen, C. The zero-free intervals for chromatic polynomials of graphs. Combin. Probab. Comput. 6 (1997), no. 4, 497–506.

[Tho00] Thomassen, C. Chromatic roots and hamiltonian paths. J. Combin. Theory B 80 (2000), 218-224.

[Tho01] Thomassen, C. Chromatic graph theory: challenges for the 21st century. World Scientific, Singapore, (2001), 183-195.

[Tra89] Traldi, L.  A dichromatic polynomial for weighted graphs and link polynomials. Proceedings of the American Mathematical Society 106 (1989), 279-286.

[Tra02] Traldi, L., Chain polynomials and Tutte polynomials. Discrete Math. 248 (2002), 279-282.

[Tra06] Traldi, L. On the colored Tutte polynomial of a graph of bounded treewidth.  Discrete Applied Mathematics 154 (2006), 1032-1036.

[TS02]      Turner, S.; Sherratt, J.A. Intercellular Adhesion and Cancer Invasion. J. Theor. Biol 216 (2002), 85-100.

[Tut47] Tutte, W.T.  A ring in graph theory.  Proc. Cambridge Philos. Soc.  43 (1947), 26-40.

[Tut53] Tutte, W.T.  A Contribution to the Theory of Chromatic Polynomials. Cand. J. Math. 6 (1953), 80-81.

[Tut67] Tutte, W.T.  On dichromatic polynomials. J.  Combin. Theory 2 (1967), 301-320.

[Tut84] Tutte, W.T. Graph Theory. Encyclopedia of Mathematics and its Applications, ed. Rota, G. C., vol. 21, Addison-Wesley, New York, 1984.

[Ver05] Vertigan, D.L.  The computational complexity of Tutte invariants for planar graphs. SIAM J. Comput. 35 (2005), no. 3, 690-712.

[VW92]  Vertigan, D.L.; Welsh, D.J.A.  The computational complexity of the Tutte plane: the bipartite case. Combin., Prob. Comp. 1 (1992), no. 2, 181-187.

[Wel93] Welsh, D.J.A.  Complexity:  Knots, colourings and counting.  London Mathematical Society Lecture Note Series.  Cambridge University Press,  New York, New York, (1993).

[Wel99] Welsh, D.J.A.  The Tutte polynomial.  Statistical physics methods in discrete probability, combinatorics, and theoretical computer science.  Random Structure Algorithms 15 (1999), no. 3-4, 210-228.

[WKS02] Wang, J.-S.; Kozan, O.; Swendsen, R.H. Sweeney and Gliozzi dynamics for simulations of Potts models in the Fortuin-Kasteleyn representation. Phys. Rev. E66 (2002), 057101-1 to 057101-4.

[WM00]  Welsh, D.J.A.; Merino, C.  The Potts model and the Tutte polynomial. Probabilistic techniques in equilibrium and nonequilibrium statistical physics.  Journal of Mathematical Physics 41 (2000), no. 3, 1127-1152.

[Woo92] Woodall, D. A zero-free interval for chromatic polynomials. Discrete Math. 101 (1992), 333-341.

[Woo02] Woodall, D. Tutte polynomials for 2-separable graphs. Discrete Math. 247 (2002) 201-213.

[Wu82]  Wu, F.-Y.  The Potts model.  Rev. Mod. Phys 54 (1982), 253-268.

[Wu84]  Wu, F.-Y. Potts model of magnetism. J. Appl. Phys. 55 (1984), 2421-2425.

[Wu88]  Wu, F.-Y. Potts model and graph theory. J. Stat. Phys. 52 (1988), 99-0112.

[Wu05]  Wu, F. Y. The random cluster model and a new integration identity.  J. Phys. A 38 (2005), 6271-6276.

[Wu06]  Wu, F.Y. New Critical Frontiers for the Potts and Percolation Models. Phys. Rev. Lett. 96 (2006), 090602.

[WW01]  Wu, F.Y.; Wang, J. Zeroes of the Jones Polynomial. Physica A 296 (2001), 483-494.

[YL52]  Yang, C.-N.; Lee, T.-D. Statistical theory of equations of state and phase transitions. I. Theory of condensation. Phys. Rev. 87 (1952), 404-409.

[Zas92] Zaslavsky, T.  Strong Tutte functions of matroids and graphs. Trans. Amer. Math. Soc. 334 (1992),  no. 1, 317-347.

\end{document}